# Splitted Levenberg-Marquardt Method for Large-Scale Sparse Problems

Nataša Krejić[*], Greta Malaspina[*], Lense Swaenen[†]

**Abstract**

We consider large-scale nonlinear least squares problems with sparse residuals, each of them depending on a small number of variables. A decoupling procedure which results in a splitting of the original problems into a sequence of independent problems of smaller sizes is proposed and analysed. The smaller size problems are modified in a way that offsets the error made by disregarding dependencies that allow us to split the original problem. The resulting method is a modification of the Levenberg-Marquardt method with smaller computational costs. Global convergence is proved as well as local linear convergence under suitable assumptions on sparsity. The method is tested on the network localization simulated problems with up to one million variables and its efficiency is demonstrated.

## 1 Introduction

The problem we consider is a nonlinear least squares problem

$$\min_{\mathbf{x}\in\mathbb{R}^N} \frac{1}{2}\sum_{j=1}^m r_j(\mathbf{x})^2 = \frac{1}{2}\min_{\mathbf{x}\in\mathbb{R}^N}\|\mathbf{R}(\mathbf{x})\|_2^2 = \min_{\mathbf{x}\in\mathbb{R}^N}\mathbf{F}(\mathbf{x}) \qquad (1)$$

[*]Department of Mathematics and Informatics, Faculty of Sciences, University of Novi Sad, Trg Dositeja Obradovića 4, 21000 Novi Sad, Serbia, Email: greta.malaspina@dmi.uns.ac.rs, natasak@uns.ac.rs.

[†]Mathware department, Sioux Technologies, Esp 130, 5633 AA Eindhoven, The Netherlands. Email: lense.swaenen@sioux.eu.

where for every $j = 1, \ldots, m$, $r_j : \mathbb{R}^N \to \mathbb{R}$ is a $C^1$ function, $\mathbf{R}(\mathbf{x}) = (r_1(\mathbf{x}), \ldots, r_m(\mathbf{x}))^\top \in \mathbb{R}^m$ is the vector of residuals, and $\mathbf{F}$ is the aggregated residual function. We are assuming a special structure of the residuals that is relevant in many applications. That is, we assume that the residuals are very sparse functions, each depending only on a small number of variables while the whole problem is large-scale, i.e. with $N$ very large. We are not assuming any particular sparsity pattern allowing many practical problems to fit the framework we consider.

Typical problems of this kind are Least Squares Network Adjustment [17], Bundle Adjustment [13], Wireless Sensor Network Localization [15], where the variables correspond to the coordinates of physical points in a region of the 2 or 3 dimensional space and residuals correspond to observations of geometrical quantities involving these points. In these cases each observation typically involves a small (often fixed) number of points, and thus each residual function involves a small number of variables. Moreover, when considering problems of large dimension with points deployed in a large region of the space, the number of observations involving each point is small with respect to the total amount of observations. That is, each variable is involved in a relatively small number of residual functions. This leads to problems that are very sparse. Given that the measurements are prone to errors or to different kinds of noise, the residuals are in general weighted with the weights being reciprocal of the measurement precision. Furthermore, a typical property of such problems is that they are nearly separable, meaning that it is possible to partition the points into subsets that are connected by a small number of observations. The dominant properties of all these problems, a very large dimension $N$ and sparsity, motivated the modification of classical Levenberg-Marquardt method presented in this paper.

The Levenberg-Marquardt (LM) method is generally used to solve the least squares problems of large dimension. The method is based on a regularization strategy for improvement of Gauss-Newton method. In each iteration of Gauss-Newton method one has to solve a linear least squares problem. The LM method further adds a regularization (damping) parameter that facilitate the direction computation. Thus in each iteration of the LM method one has to solve a system of linear equations to compute the step or step direction. The regularization parameter plays a fundamental role and its derivation is subject of many studies.

Many modifications of the classical Levenberg-Marquardt scheme have been proposed in literature to retain convergence while relaxing the assumptions on the objective function and to improve the performance of the method. In [6, 7, 18] the damping parameter is defined as a multiple of the objective function. With this choice of the parameter local superlinear or quadratic convergence is proved under a local error bound assumption for zero residual problems, while global convergence is achieved by employing a line search strategy. In [11] the authors propose an updating strategy for the parameter that, in combination with Armijo line search, ensures global convergence and $q$-quadratic local convergence under the same assumption of the previous papers. In [2] the non-zero residual case is considered and a Levenberg-Marquardt scheme is proposed that achieves local convergence with order depending on the rank of the Jacobian matrix and of a combined measure of nonlinearity and residual size around stationary points. In [4] an inexact Levenberg-Marquardt is considered and local convergence is proved under a local error bound condition. In [3] the authors propose an approximated Levenberg-Marquardt method, suitable for large-scale problems, that relies on inaccurate function values and derivatives.

The problems we are interested in are of very large dimension and sparse. Sparsity very often induces the property we call *near-separability*, i.e. it is possible to partition the variables into subsets such that a subset of residual functions depends only on a subset of variables while only a limited number of residual functions depends on variables in different subsets. This property is mainly a natural consequence of the problem origin. For example in the Network Adjustment problem physical distance of the points determines the set of points that are connected by observations. Although the complete separability of the residuals is not common, the number of residuals that depend on points from different subsets is in general rather small compared to the problem size $N$. On the other hand, for a very large $N$ solving the LM linear system in each iteration, even for very sparse problems can be costly. The particular example is the refinement of cadastral maps and in this case $N$ is prohibitively large for direct application of the LM method. For example, the Dutch Kadaster is pursuing making the cadastral map more accurate by making the map consistent with accurate field surveyor measurements [8, 10]. This application yields a non-linear least squares problem which is known as 'least squares adjustment' in the field of geography. If the entire Netherlands were to be considered for one big adjustment problem, the number of variables would be twice the number of feature points in the Netherlands, which is on the order of 1 billion variables and even considering separate parts of Netherlands still yields a very large-scale problem.

The method we propose here is designed to exploit the sparsity and near-separability in the following way. Assuming that we can split the variables in such a way that a large number of residual function depends only on a particular subset of variables, while a relatively small number of residual functions depends on variables from different subsets, the system of linear equations in LM iteration has a particular block structure. The variable subsets and corresponding residuals imply strong dominance of relatively dense diagonal blocks while the non-diagonal blocks are very sparse. This structure motivated the idea of splitting: we decompose the LM system matrix into $K$ independent systems of linear equations, determined by the diagonal blocks of the LM system. This way we get $K$ linear systems of significantly smaller dimensions and we can solve them faster, either sequentially or in parallel. Thus, one can consider this approach as a kind of inexact LM method. However, this kind of splitting might be too inaccurate given that we completely disregarded the off-diagonal blocks of the LM matrix. Therefore, we modify the $K$ independent linear systems in such a way that off-diagonal blocks are included in the right-hand sides of these independent systems. The modification of the right-hand sides of independent linear systems is based on a correction parameter proposed for the Newton method in [14] and for the distributed Newton method in [1] and attempts to minimize the difference between the full LM linear system residual and the residual of the modified method. Having an affordable search direction computed by solving the $K$ independent linear systems we can proceed in the usual way - applying a line search to get a globally convergent method. Furthermore, under a set of standard assumptions we can prove local linear convergence.

A proper splitting of the variables into suitable subsets is a key assumption for the efficiency of this method but the problems we are interested in very often offer a natural way of meaningful splitting. For example in the localization problems or network adjustment problems the geometry of points dictates meaningful subsets. In the experiments we present here one can see that a graph partitioning algorithm provides a good subset division in a cost-efficient way.

The numerical experiments presented in the paper demonstrate clearly the advantage of the proposed method with respect to the full-size LM method. We show that the splitting method successfully copes with very large dimensions, testing the problems of up to 1 million variables. Furthermore, we demonstrate that even in the case of large dimensions that can be solved by classical LM method, say around a couple of tens of thousands, the splitting method works faster. Additionally, we investigate the robustness of the splitting in the sense that we demonstrate empirically that the number of subsets plays an important role in the efficiency of the method but that there is a reasonable degree of freedom in choosing that number without affecting the performance of the method. The experiments presented here are done in a sequential way, i.e. we did not solve the independent systems of linear equations in parallel which would further enhance the proposed method. Parallel implementation will be a subject of further research.

The paper is organized as follows. In Section 2 we present the framework that we are considering. The proposed method is described in Section 3 while the theoretical analysis of the proposed method is carried out in Section 4. In Section 5 we discuss some implementation details and present numerical results.

The notation we use is the following. Mappings defined on $\mathbb{R}^N$, vectors from $\mathbb{R}^N$ and matrices with at least one dimension being $N$ are denoted by boldfaced letters $\mathbf{F}, \mathbf{R}, \mathbf{x}, \mathbf{B}, \ldots$ while their block-elements are denoted by the same letters in italics and indices so $\mathbf{x} = (x_1, \ldots, x_s)$, $\mathbf{x} \in \mathbb{R}^N$, $x_s \in \mathbb{R}^{n_s}$. The dimensions are clearly stated to avoid confusion. We use $\lambda_{\min}(\cdot)$ and $\lambda_{\max}(\cdot)$ to denote the smallest and largest eigenvalue of a matrix, respectively. The Euclidean norm is denoted by $\|\cdot\|$ for both matrices and vectors.

## 2 Nearly Separable Problems

The problem we consider is stated in (1). Denote with $\mathcal{I} = \{1, \ldots, N\}$ and with $\mathcal{J} = \{1, \ldots, m\}$. Given a partition $I_1, \ldots, I_K$ of $\mathcal{I}$ we define the corresponding partition of $\mathcal{J}$ into $E_1, \ldots, E_K$ as follows:

$$E_s = \{j \in \mathcal{J} | r_j \text{ only depends on variables in } I_s\}, \ s = 1, \ldots, K$$
$$\widehat{E} = \mathcal{J} \setminus \bigcup_{i=1}^{K} E_i. \tag{2}$$

That is, given a partition of the set of variables, each of the subsets $E_s$ contains the indices corresponding to residual functions that only involve variables in $I_s$, while $\widehat{E}$ contains the indices of residuals that involve variables belonging to different subsets $I_s$. We say that problem (1) is *separable* if there exist $K \geq 2$ and a partition $\{I_s\}_{s=1,\ldots,K}$ of $\mathcal{I}$ such that $\widehat{E} = \emptyset$, while we say that it is *nearly-separable* if there exist $K \geq 2$ and a partition $\{I_s\}_{s=1,\ldots,K}$ of $\mathcal{I}$ such that the cardinality of $\widehat{E}$ is small with respect to the cardinality of $\bigcup_{s=1}^K E_s$. The term "nearly-separable" is not defined precisely and should be understood in the same fashion as sparsity, i.e. assuming that we can identify the corresponding partitions.

The above described splitting can be interpreted as follows. Given the least squares problem in (1) we define the corresponding underlying network as the undirected graph $\mathcal{G} = (\mathcal{I}, \mathcal{E})$ where $\mathcal{I}$ and $\mathcal{E}$ denote the set of nodes and set of edges respectively. The graph $\mathcal{G}$ that has a node for each variable $x_i$, and an edge between node $i$ and node $k$ if there is a residual function $r_j$ that involves $x_i$ and $x_k$. With this definition in mind the partition of $\mathcal{I}$ and $\mathcal{J}$ that we described above corresponds to a partition of the sets of nodes and edges of the network, where $E_s$ contains the indices of edges that are between nodes in the same subset $I_s$ and $\widehat{E}$ contains the edges that connect different subsets. The problem is separable if the underlying network $\mathcal{G}$ is not connected (and the number $K$ is equal to the number of connected components of $\mathcal{G}$). The problem is nearly-separable if we can partition the set of nodes of the network in such a way that the number of edges that connect different subsets is small with respect to the number of edges that are internal to the subsets.

Given the partition $\{I_s\}_{s=1,\ldots,K}$ of the variables and the corresponding partition $\{E_s\}_{s=1,\ldots,K}$, $\widehat{E}$ of the residuals, for $s = 1, \ldots, K$ we define $x_s \in \mathbb{R}^{n_s}$ as the vector of the variables in $I_s$ where $n_s$ denotes the cardinality of $I_s$, and we introduce the following functions

$$R_s(x_s) := (r_j(\mathbf{x}))_{j \in E_s}, \qquad \rho(\mathbf{x}) := (r_j(\mathbf{x}))_{j \in \widehat{E}} \tag{3}$$
$$F_s(x_s) := \|R_s(x_s)\|_2^2 \qquad \Phi(\mathbf{x}) := \|\rho(\mathbf{x})\|_2^2$$

so that for every $s = 1, \ldots, K$, $R_s : \mathbb{R}^{n_s} \to \mathbb{R}$ is the vector of residuals involving only variables in $I_s$, while $\rho : \mathbb{R}^N \to R$ is the vector of residuals in $\widehat{E}$ and $F_s$, $\Phi$ are the corresponding local aggregated residual functions. Notice that $\sum_{s=1}^K n_s = N$. With this notation problem (1) can be rewritten as

$$\min_{\mathbf{x} \in \mathbb{R}^N} \left( \Phi(\mathbf{x}) + \sum_{s=1}^K F_s(x_s) \right) = \min_{\mathbf{x} \in \mathbb{R}^N} \left( \frac{1}{2} \|\rho(\mathbf{x})\|_2^2 + \sum_{s=1}^K \frac{1}{2} \|R_s(x_s)\|_2^2 \right) \tag{4}$$

In particular, if the problem is separable (and therefore $\widehat{E}$ is empty) $\Phi \equiv 0$ and solving problem (4) is equivalent to solving $K$ independent least squares

problems given by

$$\min_{x_s \in \mathbb{R}^{n_s}} F_s(x_s) = \min_{x_s \in \mathbb{R}^{n_s}} \frac{1}{2}\|R_s(x_s)\|_2^2 \ \text{ for } s = 1, \ldots, K. \qquad (5)$$

If the problem is not separable then in general $\Phi$ is not equal to zero and that is the case we are interested in.

## 3 LMS: the Levenberg-Marquardt method with splitting

Let $\{I_s\}_{s=1,\ldots,K}$ be a partition of $\mathcal{I}$ and $\{E_s\}_{s=1,\ldots,K}$ be the corresponding partition of $\mathcal{J}$ as defined in (2). To ease the notation we assume that we reordered the variables and the residuals functions according to the given partitions, such that for $\mathbf{x} \in \mathbb{R}^N$

$$\mathbf{x} = \begin{pmatrix} x_1 \\ \vdots \\ x_K \end{pmatrix}, \qquad \mathbf{R}(\mathbf{x}) = \begin{pmatrix} R_1(x_1) \\ \vdots \\ R_K(x_K) \\ \rho(\mathbf{x}) \end{pmatrix}$$

With this reordering, denoting with $J_{sR_j}$ the Jacobian of the partial residual vector $R_j$ defined in (3) with respect to the variables in $I_s$, and with $J_{s\rho}$ the Jacobian of the partial residual $\rho$ with respect to $x_s$, we have

$$\mathbf{J}(\mathbf{x}) = \begin{pmatrix} J_{1R_1}(x_1) & & & 0 \\ & J_{2R_2}(x_2) & & \\ & & \ddots & \\ 0 & & & J_{KR_K}(x_K) \\ J_{1\rho}(\mathbf{x}) & J_{2\rho}(\mathbf{x}) & \ldots & J_{K\rho}(\mathbf{x}) \end{pmatrix}.$$

From this structure of $\mathbf{R}$ and $\mathbf{J}$ we get the corresponding block structure of the gradient $\mathbf{g}(\mathbf{x}) = \mathbf{J}(\mathbf{x})^\top \mathbf{R}(\mathbf{x})$ and the matrix $\mathbf{J}(\mathbf{x})^\top \mathbf{J}(\mathbf{x})$:

$$\mathbf{g}(\mathbf{x})^\top = \left(g_1^\top(\mathbf{x}), g_2^\top(\mathbf{x}), \ldots, g_K^\top(\mathbf{x})\right), \qquad (6)$$

$$\mathbf{J}(\mathbf{x})^\top \mathbf{J}(\mathbf{x}) = \begin{pmatrix} H_1(\mathbf{x}) & B_{12}(\mathbf{x}) & \ldots & B_{1K}(\mathbf{x}) \\ B_{21}(\mathbf{x}) & H_2(\mathbf{x}) & \ddots & \vdots \\ \vdots & \ddots & \ddots & B_{K-1K}(\mathbf{x}) \\ B_{K1}(\mathbf{x}) & \ldots & B_{KK-1}(\mathbf{x}) & H_K(\mathbf{x}) \end{pmatrix}, \qquad (7)$$

with

$$g_s(\mathbf{x}) = J_{sR_s}(x_s)^\top R_s(x_s) + J_{s\rho_s}(\mathbf{x})^\top \rho(\mathbf{x}) \quad \text{for} \quad s = 1, \ldots, K,$$
$$H_s(\mathbf{x}) = J_{sR_s}(x_s)^\top J_{sR_s}(x_s) + J_{s\rho}(\mathbf{x})^\top J_{s\rho}(\mathbf{x}) \quad \text{for} \quad s = 1, \ldots, K, \quad (8)$$
$$B_{ij}(\mathbf{x}) = J_{i\rho}(\mathbf{x})^\top J_{j\rho}(\mathbf{x}) \quad \text{for} \quad i, j = 1, \ldots, K.$$

In the following, we denote with $\mathbf{g}^k = \mathbf{J}_k^\top \mathbf{R}_k$ the vector with $s$-th block component equal to $g_s(\mathbf{x}^k)$, with $\mathbf{H}_k = \mathbf{H}(\mathbf{x}^k)$ the block diagonal matrix with diagonal blocks given by $H_s(\mathbf{x}^k) = H_s(\mathbf{x}^k)$ for $s = 1, \ldots, K$, and with $\mathbf{B}_k = \mathbf{B}(\mathbf{x}^k)$ the block partitioned matrix with diagonal blocks equal to zero and off-diagonal blocks equal to $B_{ij}(\mathbf{x}^k)$. That is,

$$\mathbf{H}_k = \begin{pmatrix} H_1(\mathbf{x}) & & & \\ & H_2(\mathbf{x}) & & \\ & & \ddots & \\ & & & H_K(\mathbf{x}) \end{pmatrix},$$
$$\mathbf{B}_k = \begin{pmatrix} 0 & B_{12}(\mathbf{x}^k) & \ldots & B_{1K}(\mathbf{x}^k) \\ B_{21}(\mathbf{x}^k) & 0 & \ddots & \vdots \\ \vdots & \ddots & \ddots & B_{K-1K}(\mathbf{x}^k) \\ B_{K1}(\mathbf{x}^k) & \ldots & B_{KK-1}(\mathbf{x}^k) & 0 \end{pmatrix}. \quad (9)$$

The algorithm we introduce here is motivated by near-separability property and hence we state the formal assumption below.

**Assumption 1.** *There exists a constant $M > 0$ such that for all $\mathbf{x} \in \mathbb{R}^N$*

$$\|\mathbf{B}(\mathbf{x})\| \leq M \|\mathbf{J}(\mathbf{x})^\top \mathbf{J}(\mathbf{x})\|. \quad (10)$$

The assumption above is not restrictive as $\mathbf{B}(\mathbf{x})$ is a submatrix of $\mathbf{J}(\mathbf{x})^\top \mathbf{J}(\mathbf{x})$. Furthermore, the global convergence of the algorithm we propose does not depend on it in the sense that we do not use the assumption for the convergence proof. In fact the proposed algorithm works even for problems that are not nearly-separable as it can be seen as a kind of quasi-Newton method, however the efficiency of the algorithm depends on the near-separability of the problem and the value of $M$. Moreover the value of $M$ plays an important role in the analysis of local convergence and in achieving linear rate.

Consider a standard iteration of LM method for a given iteration $\mathbf{x}^k$

$$\mathbf{x}^{k+1} = \mathbf{x}^k + \mathbf{d}^k,$$

where $\mathbf{d}^k \in \mathbb{R}^N$ is the solution of

$$\left(\mathbf{J}_k^\top \mathbf{J}_k + \mu_k \mathbf{I}\right) \mathbf{d}^k = -\mathbf{J}_k^\top \mathbf{R}_k, \tag{11}$$

where $\mathbf{J}_k = \mathbf{J}(\mathbf{x}^k) \in \mathbb{R}^{m \times N}$ denotes the Jacobian matrix of $\mathbf{R}_k = \mathbf{R}(\mathbf{x}^k)$ and $\mu_k$ is a positive scalar. When $N$ is very large solving (11) at each iteration of the method may be prohibitively expensive. In the following we propose a modification of the Levenberg-Marquardt method that exploits near-separability of the problem to approximate the linear system (11) with a set of independent linear systems of smaller size.

The linear system (11) at iteration $k$ can therefore be rewritten as

$$(\mathbf{H}_k + \mu_k \mathbf{I} + \mathbf{B}_k)\mathbf{d}^k = -\mathbf{g}^k. \tag{12}$$

The matrix $\mathbf{B}_k$ depends only on the derivatives of the residual vector $\rho = (r_j(\mathbf{x}))_{j \in \widehat{E}}$. If the problem is separable then $\widehat{E} = \emptyset$ and $\mathbf{B}_k = 0$ so the coefficient matrix of (12) is block diagonal, the system can be decomposed into $K$ independent linear systems, and the solution $\mathbf{d}^k$ is a vector with the $s$-th block component equal to the solution $d_s^k$ of

$$(H_s^k + \mu_k I)d_s^k = -g_s^k. \tag{13}$$

If the problem is nearly-separable, the number of nonzero elements in $\mathbf{B}_k$ is small compared to the size $N$ of the matrix and to the number of nonzero elements in $\mathbf{H}_k$. Thus the solution of (13) may provide an approximation of the solution of the Levenberg-Marquardt direction (11), with the quality of the approximation depending on the number and magnitude of the nonzero elements in $\mathbf{B}_k$.

Given that information contained in $\mathbf{B}_k$ might be relevant and that solving $K$ systems of smaller dimension is much cheaper than solving the system of dimension $N$, we propose the following modification of the right hand side of (13), which attempts to exploit the information contained in the off-diagonal blocks, while retaining separability of the approximated linear system. The idea underlying the right-hand side correction is analogous to the one proposed in [14, 1] for systems of nonlinear equations and distributed optimization problems.

Our goal is to split the LM linear system into separable systems of smaller dimension. Starting from the LM linear system (12) and aiming at a separable system of linear equations, i.e. a system with the matrix $\mathbf{H}_k + \mu_k \mathbf{I}$ as in (13), we need to take into account the fact that $\mathbf{B}_k$ is not zero. Clearly putting $\mathbf{B}_k \mathbf{d_k}$ to the right hand side would be ideal but $\mathbf{d}_k$ is unknown. Therefore we add $\mathbf{B}_k \mathbf{g_k}$ at the right hand side of (13) as this way we maintain separability and information contained in $\mathbf{B}_k$. Intuitively, $\mathbf{B}_k \mathbf{g_k}$ is the best approximation for $\mathbf{B}_k \mathbf{d_k}$ that we have available, so we use it to get a separable system of linear equations. To compensate (at least partially) for the substitution of $\mathbf{B}_k \mathbf{d_k}$ with $\mathbf{B}_k \mathbf{g_k}$ we also add a correction factor $\beta_k$ as explained in (14) and further on.

Consider the system

$$(\mathbf{H}_k + \mu_k \mathbf{I})\mathbf{d}^k = (\beta_k \mathbf{B}_k - \mathbf{I})\mathbf{g}^k \tag{14}$$

where $\beta_k \in \mathbb{R}$ is a correction coefficient that we can choose. Once the right-hand side has been computed, since $\mathbf{H}_k + \mu_k \mathbf{I}$ is a block diagonal matrix, (14) can still be decomposed into $K$ independent linear system. That is, the solution $\mathbf{d}^k$ of (14) is given by

$$\mathbf{d}^k = \begin{pmatrix} d_1^k \\ \vdots \\ d_K^k \end{pmatrix} \quad \text{with} \quad H_s^k d_s^k = \beta_k \sum_{j=1,\ j\neq s}^{K} B_{sj}^k g_j^k - g_s^k. \tag{15}$$

The correction coefficient $\beta_k$ can be chosen freely at each iteration so far. However we will see that it is of fundamental importance for both the convergence analysis of the method and practical performance. This parameter is further specified in the algorithm we propose and discussed in detail in Subsection 4.3. Let us now give only a rough reasoning behind its introduction. With $\beta_k$ we are trying to preserve some information contained in $\mathbf{B}_k$ in a cheap way and without destroying separability. One possibility is the following choice of $\beta_k$, which ensures that the residual given by the solution of (14) with respect to the exact linear system (12) is minimized. That is,

$$\beta_k = \arg \min_{\beta \in \mathbb{R}} \|\varphi_k(\beta)\|_2^2$$

with

$$\varphi(\beta) = (\mathbf{H}_k + \mu_k \mathbf{I} + \mathbf{B}_k)(\mathbf{H}_k + \mu_k \mathbf{I})^{-1}(\beta \mathbf{B}_k - I)\mathbf{g}^k + \mathbf{g}^k. \tag{16}$$

Further details on this choice are presented later on. One can ask why we use a single coefficient $\beta_k$ in each iteration, i.e. if it would perhaps be more efficient to try to "correct" the right hand side with more than one parameter, maybe allowing a diagonal matrix $\beta_k$ in the right hand side instead of a single scalar. In fact, if we take a diagonal matrix $\beta_k \in \mathbb{R}^{N \times N}$ and plug it in the above minimization problem, then solving this problem we could recover the LM iteration exactly. However such procedure would imply the same cost for one iteration as the full LM iteration. Clearly, there are other alternatives between these two extremes - a single number and $N$ numbers, but our experience show that the choice with a single coefficient $\beta_k$ brings the best results in terms of cost benefit function. The convergence analysis presented in the next Section will further restrict the values of parameter $\beta_k$.

Consider the block diagonal matrix $\mathbf{H}_k$ defined in (7). Since the $s$-th diagonal block $H_s^k = (J_{sR_s}^k)^\top J_{sR_s}^k + (J_{s\rho_s}^k)^\top J_{s\rho_s}^k$ we have that $\mathbf{H}_k$ is symmetric and positive semi-definite, and therefore there exists a matrix $\mathbf{S}_k \in \mathbb{R}^{m \times N}$ such that $\mathbf{S}_k^\top \mathbf{S}_k = \mathbf{H}_k$. We denote with $\mathbf{C}_k$ the matrix $\mathbf{C}_k = \mathbf{J}_k - \mathbf{S}_k$.

Let us now state the algorithm. We will assume that the splitting into suitable sets is done before the iterative procedure and it is kept fixed through the process. Thus, the diagonal matrix $\mathbf{H}_k$ and off-diagonal $\mathbf{B}_k$ are already defined in each iteration.

The regularization parameter $\mu_k$ plays an important role in (14) and consequently in the resolution of our independent problem stated in (15). In the algorithm below we adopt a choice based on the same principles proposed in [11] although other options are possible. The parameter is thus computed using the values defined as

$$\hat{a}_0^k = \frac{\ell_k^2}{4} \|\mathbf{H}_k\| \|\mathbf{J}_k\| \|\mathbf{R}_k\|$$
$$\hat{a}_1^k = \frac{\ell_k^2}{4} \|\mathbf{J}_k\| \|\mathbf{R}_k\| + \ell_k \|\mathbf{H}_k\| \|\mathbf{J}_k\|^2 \|\mathbf{R}_k\| \qquad (17)$$
$$\hat{a}_2^k = \|\mathbf{H}_k\| \|\mathbf{J}_k\|^2 + \ell_k \|\mathbf{H}_k\| \|\mathbf{R}_k\| + \ell_k \|\mathbf{J}_k\|^2 \|\mathbf{R}_k\|$$
$$\hat{a}_3^k = \|\mathbf{J}_k\|^2 + \ell_k \|\mathbf{R}_k\|$$

for each iteration with $\ell_k$ specified in the algorithm below.

**Algorithm 1.**
**Parameters:**
$c \in (0, 1)$, $b \in (0, 1)$, $\eta \in [0, 1)$, $0 \leq \ell_{min} < \ell_0$.
**Iteration $k$:**
1: compute $\mu_k = \left(1 + \frac{(1+b)^2}{1-b} \max_{i=0:3} \hat{a}_i^k\right)$ with $\hat{a}_i^k$ defined in (17)

2: choose $\beta_k \in \mathbb{R}$ such that

$$\|\beta_k \mathbf{B}_k\| \leq \frac{b\mu_k}{\|\mathbf{H}_k\| + \mu_k} \qquad (18)$$

3: set $\gamma_k = \|I - \beta_k \mathbf{B}_k\|$
4: compute $\mathbf{d}^k$ by solving

$$(\mathbf{H}_k + \mu_k I)\mathbf{d}^k = (\beta_k \mathbf{B}_k - I)\mathbf{g}^k \qquad (19)$$

5: use backtracking to find the largest positive $t_k \leq \min\{1, 1/\gamma_k\}$ such that

$$\mathbf{F}(\mathbf{x}^k + t_k \mathbf{d}^k) \leq \mathbf{F}(\mathbf{x}^k) + ct_k(\mathbf{d}^k)^\top \mathbf{g}^k \qquad (20)$$

6: set $\mathbf{x}^{k+1} = \mathbf{x}^k + t_k \mathbf{d}^k$
7: **if** $t_k < \min\{1, 1/\gamma_k\}$ **then**
8: $\quad \ell_{k+1} = 2\ell_k$
9: **else**
10: $\quad$ compute $A_k = \mathbf{F}(\mathbf{x}^{k+1}) - \mathbf{F}(\mathbf{x}^k)$
11: $\quad$ compute $P_k = \mathbf{F}(\mathbf{x}^k) - \|\mathbf{R}_k + \mathbf{J}_k \mathbf{d}^k\|^2 - \mu_k \|\mathbf{d}^k\|^2$
12: $\quad$ **if** $A_k > \eta P_k$ **then**
13: $\quad\quad \ell_{k+1} = \max\{\ell_{min}, \ell_k/2\}$
14: $\quad$ **else**
15: $\quad\quad \ell_{k+1} = \ell_k$
16: $\quad$ **end if**
17: **end if**

The above algorithm has the usual structure of LM algorithm with line search with the key novelty being in (19) where we compute the search direction. The damping parameter $\mu_k$ is defined through steps 1 while the values of $\ell_k$ updated at lines 12-16 resemble trust-region approach. Roughly speaking, if the decrease is sufficient the damping parameter in the next iteration is decreased, otherwise we keep $\ell_{k+1} = \ell_k$. In fact the choice of $\mu_k$ will be crucial for the convergence proof, see Lemma 8. The correction parameter $\beta_k$ is specified in Step 2. Assuming the standard properties of the objective function, see ahead Assumption 2 and 3, one can always choose $\beta_k$ such that (18) holds. However, the value of $\beta_k$ depend on the norm of off-diagonal blocks. Thus the method essentially exploits the sparse structure we assume in this paper, see Assumption 1. The search direction $\mathbf{d}^k$ is computed in line 4. Clearly, the system (19) is in fact completely separable and solving it requires solving $K$ independent linear systems of the form (15).

The right-hand side correction vector in each system is also stated in (15). Thus, for the parameter $\beta_k$ chosen in line 2, to compute $\mathbf{d}^k$ we need to solve $K$ systems of linear equations

$$(H_s^k + \mu_k I)d_s^k = \beta_k \sum_{j=1,\ j\neq s}^{K} B_{sj}^k g_j^k - g_s^k,$$

each one of dimension $n_s$, and $\sum_{s=1}^{K} n_s = N$. These systems can be solved independently, either in parallel or sequentially, and the cost of their solving is in general significantly smaller than the cost of solving $N$ dimensional LM system of linear equations. These savings are meaningful since the direction $\mathbf{d}^k$ obtained this way is a descent direction as we will show in the convergence analysis. After that we invoke backtracking to fulfill the modified Armijo condition given in (20) and define a new iteration. Modification of the Armijo condition again depends on the norm of off-diagonal blocks as the step size is bounded above by $1/\gamma_k$. In the case of $\mathbf{B}_k = 0$, i.e. if the system is completely separable we get $\gamma_k = 1$ and the classical Armijo condition is recovered. In this case the system (19) is the classical LM system and the algorithm reduces to the classical LM with line search. On the other hand, for $\mathbf{B}^k \neq 0$ the value of $\|\mathbf{B}^k\|$ fundamentally influences the values of $\beta_k$ and $\gamma_k$ and the algorithm allows non-negligible values of $\beta_k, \gamma_k$ only if $\|\mathbf{B}^k\|$ is not loo large, i.e. if the problem has a certain level of separability

# 4 Convergence Analysis

The convergence analysis is divided in 2 parts, in Subsection 4.1. we prove that the algorithm is well defined and globally convergent under a set of standard assumptions, while the local convergence analysis is presented in Subsection 4.2. The choice of $\beta_k$ and its influence are discussed in Subsection 4.3.

## 4.1 Global Convergence

The following assumptions are regularity assumptions commonly used in LM methods

**Assumption 2.** *The vector of residuals* $\mathbf{R} : \mathbb{R}^N \to \mathbb{R}^m$ *is continuously differentiable.*

**Assumption 3.** *The Jacobian matrix $\boldsymbol{J} \in \mathbb{R}^{m \times N}$ of $\mathbf{R}$ is L-Lipschitz continuous. That is, for every $\mathbf{x}, \mathbf{y} \in \mathbb{R}^N$*

$$\|\boldsymbol{J}(\mathbf{x}) - \boldsymbol{J}(\mathbf{y})\| \leq L\|\mathbf{x} - \mathbf{y}\|.$$

For the rest of this subsection we assume that $\{\mathbf{x}^k\}$ is the sequence generated by Algorithm 1 with an arbitrary initial guess $\mathbf{x}^0 \in \mathbb{R}^N$.

The following Lemma, proved in [5], is needed for the convergence analysis.

**Lemma 1.** *[5] If Assumptions A2 and A3 hold, for every $\mathbf{x}$ and $\mathbf{y}$ in $\mathbb{R}^N$ we have*

$$\|\mathbf{R}(\mathbf{x} + \mathbf{y}) - \mathbf{R}(\mathbf{x}) - \boldsymbol{J}(\mathbf{x})\mathbf{y}\| \leq \frac{L}{2}\|\mathbf{y}\|^2. \tag{21}$$

**Lemma 2.** *Assume that $\mathbf{d}^k$ is computed as in (19) with $\beta_k$ satisfying (18), and that Assumption A2 holds. Then $\mathbf{d}^k$ is a descent direction for $\mathbf{F}$ at $\mathbf{x}^k$. Moreover the following inequalities hold*

*i)* $(\mathbf{g}^k)^\top \mathbf{d}^k \leq -\dfrac{(1-b)\|\mathbf{g}^k\|^2}{\|\mathbf{H}_k\| + \mu_k}.$

*ii)* $\mathbf{F}(\mathbf{x}^{k+1}) \leq \mathbf{F}(\mathbf{x}^k) - ct_k \dfrac{(1-b)\|\mathbf{g}^k\|^2}{\|\mathbf{H}_k\| + \mu_k}.$

*Proof.* We want to prove that $(\mathbf{g}^k)^\top \mathbf{d}^k \leq 0$ for every index $k$. By definition of $\mathbf{d}^k$ and using the fact that $\|A^{-1}\| \geq \|A\|^{-1}$ for every invertible matrix $A$, we have

$$\begin{aligned}
(\mathbf{g}^k)^\top \mathbf{d}^k &= -(\mathbf{g}^k)^\top (\mathbf{H}_k + \mu_k I)^{-1}(I - \beta_k \mathbf{B}_k)\mathbf{g}^k \\
&= -(\mathbf{g}^k)^\top (\mathbf{H}_k + \mu_k I)^{-1}\mathbf{g}^k + \beta_k (\mathbf{g}^k)^\top (\mathbf{H}_k + \mu_k I)^{-1}\mathbf{B}_k \mathbf{g}^k \\
&\leq \|\mathbf{g}^k\|^2 \left( \|\beta_k \mathbf{B}_k\| \|(\mathbf{H}_k + \mu_k I)^{-1}\| - \frac{1}{\|\mathbf{H}_k + \mu_k I\|} \right).
\end{aligned} \tag{22}$$

Since $\mu_k > 0$ and $\mathbf{H}_k$ is symmetric and positive semidefinite, we have

$$\|(\mathbf{H}_k + \mu_k I)^{-1}\| \leq \frac{1}{\lambda_{\min}(\mathbf{H}_k + \mu_k I)} \leq \frac{1}{\mu_k}$$

and

$$\frac{1}{\|\mathbf{H}_k + \mu_k I\|} = \frac{1}{\|\mathbf{H}_k\| + \mu_k}.$$

Using this two facts and the bound (18) on $\|\beta_k B_k\|$ in inequality (22), we get

$$(\mathbf{g}^k)^\top \mathbf{d}^k \leq \|\mathbf{g}^k\|^2 \left( \|\beta_k \mathbf{B}_k\| \|(\mathbf{H}_k + \mu_k I)^{-1}\| - \frac{1}{\|\mathbf{H}_k + \mu_k I\|} \right)$$

$$\leq \|\mathbf{g}^k\|^2 \left( \frac{b\mu_k}{\|\mathbf{H}_k\| + \mu_k} \frac{1}{\mu_k} - \frac{1}{\|\mathbf{H}_k\| + \mu_k} \right) = \frac{b-1}{\|\mathbf{H}_k\| + \mu_k} \|\mathbf{g}^k\|^2,$$

which is part $i)$ of the Lemma. Since $b < 1$ this also implies that $\mathbf{d}^k$ is a descent direction at iteration $k$.

By (20) we have that for every iteration index $k$

$$\mathbf{F}(\mathbf{x}^k + t_k \mathbf{d}^k) < \mathbf{F}(\mathbf{x}^k) + c t_k (\mathbf{d}^k)^\top \mathbf{g}^k.$$

Replacing $(\mathbf{g}^k)^\top \mathbf{d}^k$ with part $(i)$ of the statement we get $ii)$. □

**Remark 4.1.** Lemma 2 states that if the right-hand side correction coefficient $\beta_k$ is chosen to satisfy the condition (18), then $\mathbf{d}^k$ is a descent direction and therefore the backtracking procedure can always find a step size $t_k$ such that the Armijo condition (20) is satisfied. In particular this implies that Algorithm 1 is well defined. In Lemma 9 we will also prove that under the current assumptions the step size $t_k$ is bounded from below.

**Lemma 3.** *If Assumption A2 holds and $\mathbf{d}^k$ is the solution of (19), then for every iteration $k$ we have*

$$\|\mathbf{d}^k\| \leq \frac{\gamma_k}{\mu_k} \|\mathbf{g}^k\| \leq \frac{\gamma_k}{\mu_k} \|\boldsymbol{J}_k\| \|\mathbf{R}_k\|. \tag{23}$$

*Proof.* By definition of $\mathbf{d}^k$ and $\gamma_k$ we have

$$\|\mathbf{d}^k\| = \|(\mathbf{H}_k + \mu_k \mathbf{I})^{-1}(I - \beta_k \mathbf{B}_k)\mathbf{g}^k\|$$
$$\leq \|(\mathbf{H}_k + \mu_k \mathbf{I})^{-1}\| \|(I - \beta_k \mathbf{B}_k)\| \|\mathbf{g}^k\| \tag{24}$$
$$\leq \frac{1}{\lambda_{\min}(\mathbf{H}_k + \mu_k \mathbf{I})} \gamma_k \|\mathbf{g}^k\| \leq \frac{\gamma_k}{\mu_k} \|\mathbf{g}^k\|,$$

which is the first inequality in the thesis. The second inequality follows directly from the fact that $\mathbf{g}^k = \mathbf{J}_k^\top \mathbf{R}_k$. □

**Lemma 4.** *If Assumption A2 holds, for every $t \in [0, 1/\gamma_k]$ we have*

$$\|\mathbf{R}_k + t\boldsymbol{J}_k \mathbf{d}_k\|^2 \leq \|\mathbf{R}_k\|^2 + t(\mathbf{g}^k)^\top \mathbf{d}^k + t^2 \|\boldsymbol{J}\|^2 \|\mathbf{d}_k\|^2$$
$$- t \frac{(1-b)}{\gamma_k^2} \frac{\mu_k^2}{\|\mathbf{H}_k\| + \mu_k} \|\mathbf{d}^k\|^2 \tag{25}$$

*Proof.* By Lemma 3 we have

$$\|\mathbf{g}^k\| \geq \frac{\mu_k}{\gamma_k}\|\mathbf{d}^k\|.$$

Using this inequality in part *i)* of Lemma 2, we get

$$(\mathbf{g}^k)^\top \mathbf{d}^k \leq -\frac{1-b}{\|\mathbf{H}_k\| + \mu_k}\|\mathbf{g}^k\|^2 \leq \frac{1-b}{\gamma_k^2}\frac{\mu_k^2}{\|\mathbf{H}_k\| + \mu_k}\|\mathbf{d}^k\|^2. \qquad (26)$$

Using this inequality and the fact that $\mathbf{g}^k = \mathbf{J}_k^\top \mathbf{R}_k$ we then have

$$\|\mathbf{R}_k + t\mathbf{J}_k\mathbf{d}_k\|^2 = \|\mathbf{R}_k\|^2 + 2t(\mathbf{g}^k)^\top \mathbf{d}^k + t^2\|\mathbf{J}_k\mathbf{d}^k\|^2$$
$$\leq \|\mathbf{R}_k\|^2 + t(\mathbf{g}^k)^\top \mathbf{d}^k + t^2\|\mathbf{J}_k\|^2\|\mathbf{d}_k\|^2 - t\frac{(1-b)}{\gamma_k^2}\frac{\mu_k^2}{\|\mathbf{H}_k\| + \mu_k}\|\mathbf{d}^k\|^2 \qquad (27)$$

which is the thesis. $\square$

**Lemma 5.** *If Assumptions A2 and A3 hold, for every $t \in [0, 1/\gamma_k]$ we have*

$$\|\mathbf{R}(\mathbf{x}^k + t\mathbf{d}^k)\|^2 \leq \|\mathbf{R}_k\|^2 + t(\mathbf{g}^k)^\top \mathbf{d}^k + t\|\mathbf{d}^k\|^2\left(\frac{L^2}{4}t^3\|\mathbf{d}^k\|^2 + t\|\mathbf{J}_k\|^2\right.$$
$$\left. + Lt\|\mathbf{R}_k\| + Lt^2\|\mathbf{J}_k\|\|\mathbf{d}^k\| - \frac{1-b}{\gamma_k^2}\frac{\mu_k^2}{\|\mathbf{H}_k\| + \mu_k}\right) \qquad (28)$$

*Proof.* Let us introduce the following function from $\mathbb{R}$ to $\mathbb{R}^m$

$$\Psi(t) = \mathbf{R}(\mathbf{x}^k + t\mathbf{d}^k) - \mathbf{R}_k - t\mathbf{J}_k\mathbf{d}^k \qquad (29)$$

and let us notice that by Lemma 1 we have that $\|\Psi(t)\| \leq \frac{L}{2}t^2\|\mathbf{d}^k\|^2$. Using this bound on $\Psi$ and Lemma 4 we get

$$\|\mathbf{R}(\mathbf{x}^k + t_k\mathbf{d}^k)\|^2 = \|\Psi(t) + \mathbf{R}_k + t\mathbf{J}_k\mathbf{d}^k\|^2$$
$$\leq \|\Psi(t)\|^2 + \|\mathbf{R}_k + t\mathbf{J}_k\mathbf{d}^k\|^2 + 2\|\Psi(t)\|\|\mathbf{R}_k + t\mathbf{J}_k\mathbf{d}^k\|$$
$$\leq \frac{1}{4}L^2t^4\|\mathbf{d}^k\|^4 + Lt^2\|\mathbf{d}^k\|^2(\|\mathbf{R}_k\|^2 + t^2\|\mathbf{J}_k\|^2\|\mathbf{d}_k\|^2)^{1/2} + \|\mathbf{R}_k\|^2 \qquad (30)$$
$$+ t(\mathbf{g}^k)^\top \mathbf{d}^k + t^2\|\mathbf{J}_k\|^2\|\mathbf{d}_k\|^2 - t\frac{(1-b)}{\gamma_k^2}\frac{\mu_k^2}{\|\mathbf{H}_k\| + \mu_k}\|\mathbf{d}^k\|^2.$$

The thesis follows immediately, using the fact that

$$(\|\mathbf{R}_k\|^2 + t^2\|\mathbf{J}_k\|^2\|\mathbf{d}_k\|^2)^{1/2} \leq \|\mathbf{R}_k\| + t\|\mathbf{J}_k\|\|\mathbf{d}_k\|.$$

$\square$

**Lemma 6.** *Let us assume that Assumptions A2 and A3 hold, and let us denote with $\mu_k^*$ the largest root of the polynomial $q_k(\mu) = \sum_{j=0}^{4} a_j^k \mu^j$ with*

$$
\begin{aligned}
a_0^k &= \frac{L^2}{4}\|\mathbf{H}_k\|\|\mathbf{J}_k\|\|\mathbf{R}_k\| \\
a_1^k &= \frac{L^2}{4}\|\mathbf{J}_k\|\|\mathbf{R}_k\| + L\|\mathbf{H}_k\|\|\mathbf{J}_k\|^2\|\mathbf{R}_k\| \\
a_2^k &= \|\mathbf{H}_k\|\|\mathbf{J}_k\|^2 + L\|\mathbf{H}_k\|\|\mathbf{R}_k\| + L\|\mathbf{J}_k\|^2\|\mathbf{R}_k\| \\
a_3^k &= \|\mathbf{J}_k\|^2 + L\|\mathbf{R}_k\| \\
a_4^k &= -\frac{1-b}{\gamma_k^2}
\end{aligned}
\tag{31}
$$

*If $\mu_k > \mu_k^*$, then $t_k = \min\{1, 1/\gamma_k\}$.*

*Proof.* Using the bound to $\|\mathbf{d}_k\|$ given by Lemma 3, and the fact that $t_k \leq \min\{1, 1/\gamma_k\}$, we have

$$
\begin{aligned}
&\frac{L^2}{4}t^3\|\mathbf{d}^k\|^2 + t\|\mathbf{J}_k\|^2 + Lt\|\mathbf{R}_k\| + Lt^2\|\mathbf{J}_k\|\|\mathbf{d}^k\| - \frac{1-b}{\gamma_k^2}\frac{\mu_k^2}{\|\mathbf{H}_k\| + \mu_k} \\
&\leq t^3 \frac{L^2 \gamma_k^2}{4\mu_k^2}\|\mathbf{J}_k\|^2\|\mathbf{R}^k\|^2 + t\|\mathbf{J}_k\|^2 + Lt\|\mathbf{R}_k\| + t^2 \frac{L\gamma_k}{\mu_k}\|\mathbf{J}_k\|^2\|\mathbf{R}^k\| \\
&\quad - \frac{1-b}{\gamma_k^2}\frac{\mu_k^2}{\|\mathbf{H}_k\| + \mu_k} \\
&\leq \frac{L^2}{4\mu_k^2}\|\mathbf{J}_k\|^2\|\mathbf{R}^k\|^2 + \|\mathbf{J}_k\|^2 + L\|\mathbf{R}_k\| + \frac{L}{\mu_k}\|\mathbf{J}_k\|^2\|\mathbf{R}^k\| - \frac{1-b}{\gamma_k^2}\frac{\mu_k^2}{\|\mathbf{H}_k\| + \mu_k} \\
&= \frac{1}{\mu_k^2(\|\mathbf{H}_k\| + \mu_k)} q_k(\mu_k),
\end{aligned}
\tag{32}
$$

where $q_k$ is the polynomial with coefficients defined in (31). Together with Lemma 5, this implies that for every $t \leq \min\{1, 1/\gamma_k\}$ we have

$$
\|\mathbf{R}(\mathbf{x}^k + t\mathbf{d}^k)\|^2 \leq \|\mathbf{R}_k\|^2 + t(\mathbf{g}^k)^\top \mathbf{d}^k + t\frac{q_k(\mu_k)}{\mu_k^2(\|\mathbf{H}_k\| + \mu_k)}\|\mathbf{d}^k\|^2 \tag{33}
$$

Since $a_4^k < 0$ we have that $q(\mu) \to -\infty$ as $\mu \to +\infty$. This implies that if $\mu_k > \mu_k^*$ then $q(\mu_k) < 0$ and, by the inequality above, we have

$$
\|\mathbf{R}(\mathbf{x}^k + t\mathbf{d}^k)\|^2 \leq \|\mathbf{R}_k\|^2 + t(\mathbf{g}^k)^\top \mathbf{d}^k \tag{34}
$$

for every $t \leq \min\{1, 1/gamma_k\}$. Since $c \in (0, 1)$, this implies in particular that $\min\{1, 1/\gamma_k\}$ satisfies Armijo condition (20) and therefore $t_k = \min\{1, 1/\gamma_k\}$. □

**Lemma 7.** *If Assumptions A2 and A3 hold and at iteration $k$ we have $\ell_k \geq L$, then $t_k = \min\{1, 1/\gamma_k\}$.*

*Proof.* By the previous Lemma, in order to prove the thesis it is enough to show that when $\ell_k \geq L$ we have $\mu_k \geq \mu_k^*$ with $\mu_k^*$ largest root of the polynomial $q_k$ defined in (31). Using the Cauchy bound for the roots of polynomials [16], we have that

$$|\mu_k^*| \leq 1 + \max_{i=0:3} \frac{|a_i^k|}{|a_4^k|}. \tag{35}$$

Since $a_4^k = -\frac{1-b}{\gamma_k^2}$ and $\gamma_k \leq 1 + b$, we have that for every $k$

$$\frac{(1+b)^2}{1-b} \geq \frac{1}{|a_4^k|}.$$

Using this inequality, the definition of $\mu_k$, and the fact that $\hat{a}_i^k \geq 0$ for every $i = 0, \ldots, 3$, we have

$$\mu_k = 1 + \frac{(1+b)^2}{1-b} \max_{i=0:3} \hat{a}_i^k = 1 + \frac{(1+b)^2}{1-b} \max_{i=0:3} |\hat{a}_i^k| \geq 1 + \max_{i=0:3} \frac{|\hat{a}_i^k|}{|a_4^k|}.$$

This, together with (35) implies that $\mu_k \geq \mu_k^*$ which concludes the proof. □

**Lemma 8.** *If Assumptions A2 and A3 hold then we have $t_k = \min\{1, 1/\gamma_k\}$ for infinitely many values of $k$.*

*Proof.* By Lemma 7 we have that $t_k = \min\{1, 1/\gamma_k\}$ whenever $\ell_k \geq L$. Assume by contradiction that there exists an iteration index $\bar{k}$ such that for every $k \geq \bar{k}$ the step size $t_k$ is strictly smaller than $\min\{1, 1/\gamma_k\}$. Since in Algorithm 1 we have $\ell_{k+1} = 2\ell_k$ whenever $t_k < \min\{1, 1/\gamma_k\}$, this implies that for every $k \geq \bar{k}$ we have

$$\ell_{k+1} = 2\ell_k = 2^{k-\bar{k}}\ell_{\bar{k}}.$$

Therefore there exists $k' \geq \bar{k}$ such that $\ell_{k'} \geq L$ which implies $t_k = \min\{1, 1/\gamma_k\}$, contradicting the fact that $t_k < \min\{1, 1/\gamma_k\}$ for every $k \geq \bar{k}$. □

The above Lemma allows us to prove the first global statement below. Namely, we prove that any bounded iterative sequence has at least one accumulation point which is stationary.

**Theorem 1.** Assume that Assumptions A2, A3 hold and that $\{\mathbf{x}^k\}$ is a sequence generated by Algorithm 1 with arbitrary $\mathbf{x}^0 \in \mathbb{R}^N$. If $\{\mathbf{x}^k\}$ is bounded, then it has at least one accumulation point that is also a stationary point for $\mathbf{F}(\mathbf{x})$.

*Proof.* Since $\{\mathbf{x}^k\} \subset \mathbb{R}^n$ is bounded and by Lemma 8 the sequence of step sizes $\{t_k\}$ takes value $\min\{1, 1/\gamma_k\}$ infinitely many times, we can take a subsequence $\{\mathbf{x}^{k_j}\} \subset \{\mathbf{x}^k\}$ such that $t_{k_j} = \min\{1, 1/\gamma_{k_j}\}$ for every $j$ and that $\mathbf{x}^{k_j}$ converges to $\bar{\mathbf{x}}$ as $j$ tends to infinity. By Lemma 2 we have

$$c(1-b) \sum_{j=0}^{\infty} \frac{1}{\min\{1, 1/\gamma_{k_j}\}} \frac{\|\mathbf{g}_{k_j}\|^2}{\|\mathbf{H}_{k_j}\| + \mu_{k_j}} \leq \mathbf{F}_0 < \infty$$

which implies that

$$\lim_{j \to +\infty} \frac{1}{\min\{1, 1/\gamma_{k_j}\}} \frac{\|\mathbf{g}_{k_j}\|^2}{\|\mathbf{H}_{k_j}\| + \mu_{k_j}} = 0. \tag{36}$$

By definition of $\gamma_k$ and (18) we have that $\min\{1, 1/\gamma_k\} \leq (1+b)$. Since $\{\mathbf{x}_{k_j}\}$ is a compact subset of $\mathbb{R}^N$, and $\mathbf{R}(\mathbf{x})$ is twice continuously differentiable, we have that the sequences $\|\mathbf{H}_{k_j}\|$, $\|\mathbf{R}_{k_j}\|$, and $\|\mathbf{J}_{k_j}\|$ are bounded from above, which by definition of $\mu_k$ imply that $\|\mathbf{H}_{k_j}\| + \mu_{k_j}$ is also bounded from above. This, together with (36) implies that $\|\mathbf{g}_{k_j}\|$ vanishes as $j$ tends to infinity and therefore $\bar{\mathbf{x}}$ is a stationary point of $\mathbf{F}(\mathbf{x})$. □

**Lemma 9.** *If Assumptions A2 and A3 hold and $\ell_{min} > 0$ then for every index $k$ we have*

$$t_k \geq \min\left\{\frac{\ell_{min}^8}{4L^8}, \frac{1}{1+b}\right\}.$$

*Proof.* From inequality (32), since $t \leq \min\{1, 1/\gamma_k\}$, we have

$$\frac{L^2}{4}t^3\|\mathbf{d}^k\|^2 + t\|\mathbf{J}_k\|^2 + Lt\|\mathbf{R}_k\| + Lt^2\|\mathbf{J}_k\|\|\mathbf{d}^k\| - \frac{1-b}{\gamma_k^2}\frac{\mu_k^2}{\|\mathbf{H}_k\| + \mu_k}$$
$$\leq \frac{1}{\mu_k^2(\|\mathbf{H}_k\| + \mu_k)}\left(\left(q_k(\mu_k) + \frac{1-b}{\gamma_k^2}\mu_k^4\right)t - \frac{1-b}{\gamma_k^2}\mu_k^4\right), \tag{37}$$

where $q_k$ is defined in (31). Using the inequality above together with Lemma 7 we have

$$\|\mathbf{R}(\mathbf{x}^k + t\mathbf{d}^k)\|^2 \leq \|\mathbf{R}_k\|^2 + t(\mathbf{g}^k)^\top \mathbf{d}^k$$
$$+ t\|\mathbf{d}^k\|^2 \frac{1}{\mu_k^2(\|\mathbf{H}_k\| + \mu_k)} \left(\left(q_k(\mu_k) + \frac{1-b}{\gamma_k^2}\mu_k^4\right) t - \frac{1-b}{\gamma_k^2}\mu_k^4\right) \quad (38)$$

Let us define

$$\bar{t}_k := \frac{\mu_k^4}{\frac{\gamma_k^2}{1-b} q_k(\mu_k) + \mu_k^4}.$$

We can easily see that if $t \leq \bar{t}_k$ then the term between parentheses of the previous inequality is non-positive and therefore

$$\|\mathbf{R}(\mathbf{x}^k + t\mathbf{d}^k)\|^2 \leq \|\mathbf{R}_k\|^2 + t(\mathbf{g}^k)^\top \mathbf{d}^k. \quad (39)$$

Since in Algorithm 1 we take $c \in (0,1)$ this implies that if $t_k > \bar{t}_k$ Armijo condition (20) holds.
By Lemma 8 we have that if $\ell_k \geq L$ then $t_k = \min\{1, 1/\gamma_k\} \geq 1/(b+1)$.
Let us consider the case when $\ell_k < L$, which also implies $\mu_k \leq \bar{\mu}_k$ with

$$\bar{\mu}_k = 1 + \max_{i=0:3} \frac{|a_i^k|}{|a_4^k|}.$$

Using the definition of $\bar{\mu}_k$ and the fact that $\bar{\mu}_k \geq 1$ and $\mu_k \leq \bar{\mu}_k$, we have

$$\frac{\gamma_k^2}{1-b} q_k(\mu_k) + \mu_k^4 = -\mu_k^4 + \frac{a_3^k}{|a_4^k|}\mu_k^3 + \frac{a_2^k}{|a_4^k|}\mu_k^2 + \frac{a_1^k}{|a_4^k|}\mu_k + \frac{a_0^k}{|a_4^k|} + \mu_k^4$$
$$\leq \frac{a_3^k}{|a_4^k|}\bar{\mu}_k^3 + \frac{a_2^k}{|a_4^k|}\bar{\mu}_k^2 + \frac{a_1^k}{|a_4^k|}\bar{\mu}_k + \frac{a_0^k}{|a_4^k|} \quad (40)$$
$$\leq 4\left(1 + \max_{i=0:3} \frac{|a_i^k|}{|a_4^k|}\right) \bar{\mu}_k^3 = 4\bar{\mu}_k^4.$$

Since we are considering the case $\ell_k < L$, we have $\hat{a}_i^k \leq \frac{\ell_k^2}{L^2}|a_i^k|$ for every $i = 0, \ldots, 3$. Moreover, $\frac{1}{|a_4^k|} = \frac{\gamma_k^2}{1-b} \leq \frac{(1+b)^2}{1-b}$. This implies

$$\mu_k = 1 + \frac{(1+b)^2}{1-b} \max_{i=0:3} \hat{a}_i^k \geq 1 + \frac{(1+b)^2}{1-b} \frac{\ell_k^2}{L^2} \max_{i=0:3} |a_i^k|$$
$$\geq 1 + \frac{\ell_k^2}{L^2} \max_{i=0:3} \frac{|a_i^k|}{|a_4^k|} \geq \frac{\ell_k^2}{L^2}\left(1 + \max_{i=0:3} \frac{|a_i^k|}{|a_4^k|}\right) = \frac{\ell_k^2}{L^2}\bar{\mu}_k.$$

Using this inequality and (40) in the definition of $\bar{t}_k$ we get

$$\bar{t}_k = \frac{\mu_k^4}{\frac{\gamma_k^2}{1-b}q_k(\mu_k) + \mu_k^4} \geq \frac{\mu_k^4}{4\bar{\mu}_k^4} \geq \left(\frac{\ell_k^2}{L^2}\bar{\mu}_k\right)^4 \frac{1}{4\bar{\mu}_k^4} \geq \frac{\ell_{min}^8}{4L^8} \quad (41)$$

which gives us the thesis.

$\square$

Finally, we can state the global convergence results.

**Theorem 2.** If Assumptions A2 and A3 hold and $\ell_{min} > 0$ then every accumulation point of the sequence $\{\mathbf{x}^k\}$ is a stationary point of $\mathbf{F}(\mathbf{x})$.

*Proof.* Let $\bar{\mathbf{x}}$ be an accumulation point of $\{\mathbf{x}^k\}$ and let $\{\mathbf{x}^{k_j}\}$ be a subsequence converging to $\bar{\mathbf{x}}$. By Lemma 2 we have

$$c(1-b)\sum_{j=0}^{\infty} t_{k_j} \frac{\|\mathbf{g}_{k_j}\|^2}{\|\mathbf{H}_{k_j}\| + \mu_{k_j}} \leq \mathbf{F}_0 < \infty$$

and therefore that

$$\lim_{j \to =\infty} t_{k_j} \frac{\|\mathbf{g}_{k_j}\|^2}{\|\mathbf{H}_{k_j}\| + \mu_{k_j}} = 0.$$

By Lemma 9 the sequence $\{t_k\}$ is bounded from below while by continuity of $\mathbf{J}(\mathbf{x}), \mathbf{R}(\mathbf{x}), \mathbf{H}(\mathbf{x})$, and of the norm 2, we have that $\|\mathbf{H}_{k_j}\| + \mu_{k_j}$ is bounded from above. This implies

$$0 = \lim_{j \to +\infty} \|\mathbf{g}_{k_j}\| = \|\mathbf{g}(\bar{\mathbf{x}})\|$$

and thus $\bar{\mathbf{x}}$ is a stationary point of $\mathbf{F}(\mathbf{x})$.

$\square$

## 4.2 Local Convergence

Let us now analyze the local convergence. We are going to show that the LMS method generates a linearly converging sequence under a set of suitable assumptions. Notice that the assumptions we use are standard, see [2] plus the sparsity assumption that we already stated. Let $S$ denote the set of all stationary points of $\|\mathbf{R}(\mathbf{x})\|$, namely $S = \{\mathbf{x} \in \mathbb{R}^N | \mathbf{J}(\mathbf{x})^\top \mathbf{R}(\mathbf{x}) = 0\}$. Consider a stationary point $\mathbf{x}^* \in S$ and a ball $B(\mathbf{x}^*, r)$ with radius $r > 0$ around it. In the rest of the section we make the following assumptions, see [2].

**Assumption 4.** *There exists $\omega > 0$ such that for every $\mathbf{x} \in B(\mathbf{x}^*, r)$*

$$\omega \operatorname{dist}(\mathbf{x}, S) \leq \|\boldsymbol{J}(\mathbf{x})^\top \mathbf{R}(\mathbf{x})\|$$

**Assumption 5.** *There exists $\sigma > 0$ such that for every $\mathbf{x} \in B(\mathbf{x}^*, r)$ and every $\bar{\mathbf{z}} \in B(\mathbf{x}^*, r) \cap S$*

$$\|(\boldsymbol{J}(\mathbf{x}) - \boldsymbol{J}(\bar{\mathbf{z}}))^\top \mathbf{R}(\bar{\mathbf{z}})\| \leq \sigma \|\mathbf{x} - \bar{\mathbf{z}}\|.$$

From now on we denote with $\rho_k$ the relative residual of the linear system (11) by the approximate solution $\mathbf{d}^k$. That is

$$\|\mathbf{g}^k + (\mathbf{J}_k^\top \mathbf{J}_k + \mu_k \mathbf{I})\mathbf{d}^k\| \leq \rho_k \|\mathbf{g}^k\|. \tag{42}$$

This residual is already considered in (16), where we briefly mentioned that we will determine $\beta_k$ such that this residual is minimized. Further details on this choice are presented in Subsection 4.3 and in this part we will keep $\rho_k$ without further specification, i.e., assuming only that it is small enough that local convergence requirements can be fulfilled. Clearly, for the completely separable problems, i.e. $\mathbf{B}_k = 0$ we get $\rho_k = 0$ and hence the value of $\rho_k$ depends on $M$ stated in Assumption 1 - if $M$ is small enough, i.e., if the problem is nearly-separable to the sufficient degree it is reasonable to expect that the values of $\rho_k$ will be small enough with a suitable choice of $\beta_k$.

The inequalities in the Lemma below are direct consequences of Assumption 3, their proofs can be found in [2].

**Lemma 10.** *Let Assumption 2 hold. There exist positive constants $L_2, L_3$ and $L_4$ such that for every $\mathbf{x}, \mathbf{y} \in D_r$, $\bar{\mathbf{z}} \in B(\mathbf{x}^*, r) \cap S$ the inequalities below hold:*

$$\|\mathbf{R}(\mathbf{x}) - \mathbf{R}(\mathbf{y}) - \boldsymbol{J}(\mathbf{y})(\mathbf{x} - \mathbf{y})\| \leq \frac{L}{2} \|\mathbf{x} - \mathbf{y}\|^2 \tag{43}$$

$$\|\mathbf{R}(\mathbf{x}) - \mathbf{R}(\mathbf{y})\| \leq L_2 \|\mathbf{x} - \mathbf{y}\| \tag{44}$$

$$\|\boldsymbol{J}(\mathbf{x})^\top \mathbf{R}(\mathbf{x}) - \boldsymbol{J}(\mathbf{y})^\top \mathbf{R}(\mathbf{y})\| \leq L_3 \|\mathbf{x} - \mathbf{y}\| \tag{45}$$

$$\begin{aligned} \|\boldsymbol{J}(\mathbf{x})^\top \mathbf{R}(\mathbf{x}) - \boldsymbol{J}(\mathbf{y})^\top \mathbf{R}(\mathbf{y}) - \boldsymbol{J}(\mathbf{x})^\top \boldsymbol{J}(\mathbf{x})(\mathbf{x} - \mathbf{y})\| \\ \leq L_4 \|\mathbf{x} - \mathbf{y}\|^2 + \|(\boldsymbol{J}(\mathbf{x}) - \boldsymbol{J}(\mathbf{y}))^\top \mathbf{R}(\mathbf{y})\| \end{aligned} \tag{46}$$

$$\begin{aligned} \|(\boldsymbol{J}(\mathbf{x}) - \boldsymbol{J}(\mathbf{y}))^\top \mathbf{R}(\mathbf{y})\| \leq LL_2 \|\mathbf{x} - \bar{\mathbf{z}}\| \|\mathbf{y} - \bar{\mathbf{x}}\| + LL_2 \|\mathbf{y} - \bar{\mathbf{z}}\|^2 \\ + \|\boldsymbol{J}(\mathbf{x})^\top \mathbf{R}(\bar{\mathbf{z}})\| \|\boldsymbol{J}(\mathbf{y})^\top \mathbf{R}(\bar{\mathbf{z}})\| \end{aligned} \tag{47}$$

From now on, given any point $\mathbf{x} \in \mathbb{R}^N$, we denote with $\bar{\mathbf{x}}$ the point in $S$ that realizes $\|\mathbf{x} - \bar{\mathbf{x}}\| = \text{dist}(\mathbf{x}, S)$.

**Lemma 11.** *There exists $r > 0$ and $c_1 > 0$ such that if $\mathbf{x}^k \in B(\mathbf{x}^*, r)$ then $\|\mathbf{d}^k\| \leq c_1 \text{dist}(\mathbf{x}^k, S)$.*

*Proof.* Let us define $\mathbf{H}^* = \mathbf{H}(\mathbf{x}^*)$. Consider the eigendecomposition of $\mathbf{H}^* = \mathbf{S}_*^\top \mathbf{S}_* = \mathbf{Q}_* \mathbf{\Lambda}_* \mathbf{Q}_*^\top$ where $\mathbf{\Lambda}_*$ is a diagonal matrix containing the ordered eigenvalues of $\mathbf{S}_*^\top \mathbf{S}_*$ and $\mathbf{Q}_*$ is the matrix containing the orthogonal eigenvectors corresponding to the eigenvalues in $\mathbf{\Lambda}_*$. Denoting with $p$ the rank of $\mathbf{S}_*^\top \mathbf{S}_*$ we have that

$$\mathbf{\Lambda}_* = \begin{pmatrix} \Lambda_{*1} & 0 \\ 0 & 0 \end{pmatrix}$$

with $\Lambda_{*1} = \text{diag}(\lambda_1^*, \ldots, \lambda_p^*)$ with $\lambda_1 \geq \lambda_2 \geq \cdots \geq \lambda_p > 0$. Consider the eigendecomposition of $\mathbf{H}_k = \mathbf{Q}_k \mathbf{\Lambda}_k \mathbf{Q}_k^\top$ and consider the partition of $\mathbf{Q}_k$ and $\mathbf{\Lambda}_k$ corresponding to the partition of $\mathbf{\Lambda}_*$:

$$\mathbf{Q}_k = (Q_{k1}, Q_{k2}), \quad \mathbf{\Lambda}_k = \begin{pmatrix} \Lambda_{k1} & 0 \\ 0 & \Lambda_{k2} \end{pmatrix}$$

with $\Lambda_{k1} = \text{diag}(\lambda_1^k, \ldots, \lambda_p^k) \in \mathbb{R}^{p \times p}$, $\Lambda_{k2} = \text{diag}(\lambda_{p+1}^k, \ldots, \lambda_N^k) \in \mathbb{R}^{(N-p) \times (N-p)}$, $Q_{k1} \in \mathbb{R}^{n \times p}$, $Q_{k2} \in \mathbb{R}^{n \times (N-p)}$ and $\lambda_1^k \geq \cdots \geq \lambda_p^k$. Since $\mathbf{R}$ is continuously differentiable on $B(\mathbf{x}^*, r)$ the entries of $\mathbf{S}(\mathbf{x})^\top \mathbf{S}(\mathbf{x})$ are continuous functions of $\mathbf{x}$ and thus the eigenvalues of $\mathbf{S}(\mathbf{x})^\top \mathbf{S}(\mathbf{x})$ are continuous function of $\mathbf{x}$. Therefore, for $r$ small enough we have $\lambda_i^k \geq \frac{1}{2} \lambda_p^*$ for every $i = 1, \ldots, p$. Moreover, since $\mathbf{Q}_k$ is an orthogonal matrix, we have that $\|\mathbf{d}^k\|^2 = \|Q_{k1}^\top \mathbf{d}^k\|^2 + \|Q_{k2}^\top \mathbf{d}^k\|^2$. For $i = 1, 2$, by definition of $\mathbf{d}^k$, we have

$$\mathbf{Q}_{ki}(\beta_k \mathbf{B}_k - \mathbf{I})\mathbf{g}^k = \mathbf{Q}_{ki}(\mathbf{H}_k + \mu \mathbf{I})\mathbf{d}^k = (\mathbf{\Lambda}_{ki} + \mu \mathbf{I})\mathbf{Q}_{ki}^\top \mathbf{d}^k \qquad (48)$$

so that

$$\|\mathbf{Q}_{ki}^\top \mathbf{d}^k\| = \|(\mathbf{\Lambda}_{ki} + \mu \mathbf{I})^{-1}(\beta_k \mathbf{B}_k - \mathbf{I})\mathbf{g}^k\|. \qquad (49)$$

By definition of $\gamma_k$, inequality (45), and the fact that $\lambda_p^k \geq \frac{1}{2} \lambda_p^*$ we have

$$\|Q_{k1}^\top \mathbf{d}^k\| \leq \|(\mathbf{\Lambda}_{k1} + \mu \mathbf{I})^{-1}\| \|(\beta_k \mathbf{B}_k - \mathbf{I})\| \|\mathbf{g}^k\| \leq \frac{\gamma_k}{\lambda_p^k + \mu_k} \|\mathbf{g}^k\|$$
$$\leq \frac{2\gamma_k L_3}{\lambda_p^*} \|\mathbf{x}^k - \bar{\mathbf{x}}^k\| \leq \frac{2\gamma_k L_3}{\lambda_p^*} \text{dist}(\mathbf{x}^k, S) \qquad (50)$$

and analogously,

$$\|Q_{k2}^\top \mathbf{d}^k\| \leq \frac{\gamma_k}{\mu_k}\|\mathbf{g}^k\| \leq \frac{\gamma_k L_3}{\mu_k}\operatorname{dist}(\mathbf{x}^k, S).$$

Therefore the thesis holds with

$$c_1 = \gamma L_3 \left(\frac{4}{(\lambda_p^*)^2} + \frac{1}{\mu_{\min}^2}\right)^{1/2}$$

for $\gamma = 1 + b \geq \gamma_k$. $\square$

**Lemma 12.** *If $\mathbf{x}^k, \mathbf{x}^{k+1} \in B(\mathbf{x}^*, r/2)$ then*

$$\begin{aligned}\omega \operatorname{dist}(\mathbf{x}^{k+1}, S) &\leq (L_4 c_1^2 + L L_2(2+c_1)(1+c_1))\|\mathbf{x}^k - \bar{\mathbf{x}}^k\|^2 \\ &+ (\mu_k c_1 + \rho_{\max} L_3)\|\mathbf{x}^k - \bar{\mathbf{x}}^k\| + \|\mathbf{J}_k^\top \mathbf{R}(\bar{\mathbf{x}}^k)\| + \|(\mathbf{J}_{k+1})^\top \mathbf{R}(\bar{\mathbf{x}}^k)\|\end{aligned} \quad (51)$$

*where $\rho_{max} = \max_k \rho_k$ with $\rho_k$ defined in (42) and $\bar{\mathbf{x}}^k$ is a point in $S$ such that $\operatorname{dist}(\mathbf{x}^k, S) = \|\mathbf{x}^k - \bar{\mathbf{x}}^k\|$.*

*Proof.* Since

$$\begin{aligned}\omega \operatorname{dist}(\mathbf{x}^{k+1}, S) &\leq \|\mathbf{g}_{k+1}\| \leq \|\mathbf{g}_k + \mathbf{J}_k^\top \mathbf{J}_k \mathbf{d}^k\| + \|\mathbf{g}_{k+1} - \mathbf{g}_k - (\mathbf{J}_k)^\top \mathbf{J}_k(\mathbf{x}^{k+1} - \mathbf{x}^k)\| \\ &\leq L_4 \|\mathbf{d}^k\|^2 + \|(\mathbf{J}_k - \mathbf{J}_{k+1})^\top R_{k+1}\| + \mu_k \|\mathbf{d}^k\| \\ &\quad + \|\mathbf{g}^k + (\mathbf{J}_k^\top \mathbf{J}_k + \mu_k \mathbf{I})\mathbf{d}^k\|.\end{aligned} \quad (52)$$

we have

$$\begin{aligned}\|(\mathbf{J}_k - \mathbf{J}_{k+1})^\top R_{k+1}\| &\leq L L_2(\|\mathbf{x}^k - \bar{\mathbf{x}}^k\| + \|\mathbf{x}^{k+1} - \bar{\mathbf{x}}^k\|)\|\mathbf{x}^{k+1} - \bar{\mathbf{x}}^k\| \\ &\quad + \|\mathbf{J}_k^\top \mathbf{R}(\bar{\mathbf{x}}^k)\| + \|(\mathbf{J}_{k+1})^\top \mathbf{R}(\bar{\mathbf{x}}^k)\| \\ &\leq L L_2(1+c_1)(2+c_1)\|\mathbf{x}^{k+1} - \bar{\mathbf{x}}^k\|^2 \\ &\quad + \|\mathbf{J}_k^\top \mathbf{R}(\bar{\mathbf{x}}^k)\| + \|(\mathbf{J}_{k+1})^\top \mathbf{R}(\bar{\mathbf{x}}^k)\|.\end{aligned} \quad (53)$$

By definition of $\rho_k$ there holds

$$\|\mathbf{g}^k + (\mathbf{J}_k^\top \mathbf{J}_k + \mu_k \mathbf{I})\mathbf{d}^k\| \leq \rho_{\max}\|\mathbf{g}^k\| \leq \rho_{\max} L_3 \|\mathbf{x}^k - \bar{\mathbf{x}}^k\| \quad (54)$$

Replacing these two inequalities in (52) and using Lemma 11 we get the thesis. $\square$

**Lemma 13.** *Assume that there exists $\eta \in (0,1)$ such that $\eta\omega > c_1\mu_{max} + \rho_{\max}L_3 + (2+c_1)\sigma$,*

$$\varepsilon = \frac{\eta\omega - (c_1\mu_{max} + \rho_{\max}L_3 + (2+c_1)\sigma)}{L_4c_1^2 + LL_2(1+c_1)(2+c_2)}.$$

*If $\mathbf{x}^k, \mathbf{x}^{k+1} \in B(\mathbf{x}^*, r/2)$ and $\mathrm{dist}(\mathbf{x}^k, S) \leq \varepsilon$ then*

$$\mathrm{dist}(\mathbf{x}^{k+1}, S) \leq \eta\, \mathrm{dist}(\mathbf{x}^k, S).$$

*Proof.* By Assumption 5 and Lemma 11

$$\|\mathbf{J}_k^\top \mathbf{R}(\bar{\mathbf{x}}^k)\| \leq \|(\mathbf{J}_k - \mathbf{J}(\bar{\mathbf{x}}^k))^\top \mathbf{R}(\bar{\mathbf{x}}^k)\| \leq \sigma\|\mathbf{x}^k - \bar{\mathbf{x}}^k\|$$

and

$$\|(\mathbf{J}_{k+1})^\top \mathbf{R}(\bar{\mathbf{x}}^k)\| \leq \|(\mathbf{J}_{k+1} - \mathbf{J}(\bar{\mathbf{x}}^k))^\top \mathbf{R}(\bar{\mathbf{x}}^k)\| \leq \sigma(1+c_1)\|\mathbf{x}^k - \bar{\mathbf{x}}^k\|$$

Therefore, from Lemma 12, since we are assuming $\mathrm{dist}(\mathbf{x}^k, S) \leq \varepsilon$, there follows

$$\omega\, \mathrm{dist}(\mathbf{x}^{k+1}, S) \leq (L_4c_1^2 + LL_2(2+c_1)(1+c_1))\|\mathbf{x}^k - \bar{\mathbf{x}}^k\|^2$$
$$+ (\mu_k c_1 + \rho_{\max}L_3 + (2+c_1)\sigma)\|\mathbf{x}^k - \bar{\mathbf{x}}^k\|$$
$$\leq \left((L_4c_1^2 + LL_2(2+c_1)(1+c_1))\varepsilon + (\mu_k c_1 + \rho_{\max}L_3 + (2+c_1)\sigma)\right)\|\mathbf{x}^k - \bar{\mathbf{x}}^k\|$$
(55)

and we get the thesis by definition of $\varepsilon$. □

The above Lemmas allow us to prove the local linear convergence.

**Theorem 3.** *Assume that Assumptions 2-5 hold and that there exists $\eta \in (0,1)$ such that $\eta\omega > \mu_{\max}c_1 + L_3\rho_{\max} + (2+c_1)\sigma$ and let us define*

$$\varepsilon = \min\left\{\frac{\eta\omega - (c_1\mu_{max} + \rho_{\max}L_3 + (2+c_1)\sigma)}{L_4c_1^2 + LL_2(1+c_1)(2+c_2)}, \frac{1}{2}\frac{r(1-\eta)}{1+c_1-\eta}\right\}.$$

*If $\mathbf{x}^0 \in B(\mathbf{x}^*, \varepsilon)$ then $\mathrm{dist}(\mathbf{x}^k, S) \to 0$ linearly and $\mathbf{x}^k \to \bar{\mathbf{x}} \in S \cap B(\mathbf{x}^*, r/2)$.*

*Proof.* We prove by induction on $k$ that $\mathbf{x}^k \in B(\mathbf{x}^*, r/2)$ for every $k$.
For $k=1$, by Lemma 11 and the definition of $\varepsilon$, we have

$$\|\mathbf{x}^1 - \mathbf{x}^*\| \leq \|\mathbf{x}^1 - \mathbf{x}^0\| + \|\mathbf{x}^0 - \mathbf{x}^*\| \leq \|\mathbf{d}_o\| + \varepsilon \leq \varepsilon(1+c_1) \leq \frac{r}{2}.$$

Given $k \geq 1$, assume that for every $j = 1, \ldots, k-1$ there holds $\text{dist}(\mathbf{x}^j, S) \leq \varepsilon$ and $\mathbf{x}^j \in B(x^*, r/2)$. Then we have

$$\|\mathbf{x}^{k+1} - \mathbf{x}^*\| \leq \|\mathbf{x}^1 - \mathbf{x}^*\| + \sum_{j=1}^{k} \|\mathbf{d}^j\|$$

and the fact that the right-hand side is smaller than $r/2$ follows again from Lemma 11 and the definition of $\varepsilon$. So, $\mathbf{x}^k \in B(\mathbf{x}^*, r/2)$ for every $k$ and to prove the first part of the thesis it is enough to apply Lemma 13.

Therefore, if there exists $\bar{\mathbf{x}} = \lim \mathbf{x}^k$, then the limit has to belong to $S \cap B(x^*, r/2)$, and to prove the second part of the thesis we only need to prove that such a limit exists. For every index $k$ we have that $\|\mathbf{d}^k\| \leq c_1 \varepsilon \eta^k$, so given any two indices $l, q$, we have

$$\|\mathbf{x}^l - \mathbf{x}^q\| \leq \sum_{j=q}^{l} \|\mathbf{d}_j\| \leq \sum_{j=q}^{\infty} \|\mathbf{d}_j\| \leq c_1 \varepsilon \sum_{j=q}^{\infty} \eta^j$$

and $\{\mathbf{x}^k\}$ is a Cauchy sequence in $\mathbb{R}^N$. So, it is convergent.

□

**Remark 4.2.** We notice that the condition

$$\eta \omega > \mu_{\max} c_1 + L_3 \rho_{\max} + (2 + c_1) \sigma$$

in Theorem 3 is analogous to the condition used to prove local linear convergence in [2], namely $\eta \omega > (2 + c_1) \sigma$. The two additional terms in the condition in Theorem 3 are a consequence of the main differences between Algorithm 1 and the method considered in [2]. In particular $\mu_{\max} c_1$ depends on the different choice of $\mu_k$, while $L_3 \rho_{\max}$ arise from the fact that at each iteration the Levenberg Marquardt system is solved inexactly.

We also notice that, recalling the definition of $c_1$ in Lemma 11, the condition above implies

$$\sigma \leq \frac{1}{2 + c_1}(\eta \omega - c_1 \mu_{\max} - L_3 \rho_{\max}) \leq \frac{\eta \omega}{c_1} < \frac{\eta \omega}{2 \gamma L_3} \lambda_p^* < \lambda_p^*,$$

which in turn is analogous to the condition $\sigma < \lambda_n^*$ used for the convergence analysis of classical Levenberg-Marquardt method in the case of problems with nonsingular Jacobian and nonzero residual at the solution [5].

## 4.3 Choice of $\beta_k$

The choice of $\beta$ is mentioned several times as a crucial ingredient of the algorithm we consider. Recall that the role of $\beta_k$ is to compensate, if possible, information that we disregarded by splitting the original LM system into $k$ separable systems in a computationally efficient way. Furthermore, due to condition (18) $\beta_k$ can have non-negligible value only if $\|\mathbf{B}_k\|$ is not too large, i.e., if the problem is sparse enough and that is enough for global convergence. To obtain local linear convergence we need to make the residual small enough, recall (42). An intuitive approach is to determine $\beta_k$ such that the residual given by the solution of (14) with respect to the exact linear system (12) is minimized. That is, we have

$$\beta_k = \arg\min_{\beta \in \mathbb{R}} \|\varphi_k(\beta)\|_2^2 \tag{56}$$

with

$$\varphi_k(\beta) = (\mathbf{H}_k + \mathbf{B}_k + \mu_k \mathbf{I})(\mathbf{H}_k + \mu_k \mathbf{I})^{-1}(\beta \mathbf{B}_k - I)\mathbf{g}^k + \mathbf{g}^k. \tag{57}$$

Defining $\mathbf{u} = \mathbf{B}_k \mathbf{g}^k$, $\mathbf{v} = \mathbf{B}_k(\mathbf{H}_k + \mu_k \mathbf{I})^{-1}\mathbf{B}_k \mathbf{g}^k$, $\mathbf{w} = \mathbf{B}_k(\mathbf{H}_k + \mu_k \mathbf{I})^{-1}\mathbf{g}^k$, we have

$$\|\varphi_k(\beta)\|^2 = \beta^2 \|\mathbf{u} + \mathbf{w}\|^2 - 2\beta (\mathbf{u} + \mathbf{v})^\top \mathbf{w} + \|\mathbf{w}\|^2$$

and the solution of (56) is given by

$$\beta_k = \frac{(\mathbf{u} + \mathbf{v})^\top \mathbf{w}}{\|\mathbf{u} + \mathbf{v}\|^2}. \tag{58}$$

Let us now consider the actual computation of $\beta_k$. To compute the vectors $\mathbf{u}$, $\mathbf{v}$, $\mathbf{w}$ we first compute $\mathbf{u} = \mathbf{B}_k \mathbf{g}^k$ directly, then we find $\widehat{\mathbf{v}}$, $\widehat{\mathbf{w}}$ such that $(\mathbf{H}_k + \mu_k \mathbf{I})\widehat{\mathbf{v}} = \mathbf{u}$ and $(\mathbf{H}_k + \mu_k \mathbf{I})\widehat{\mathbf{w}} = \mathbf{g}^k$, and finally we take $\mathbf{v} = \mathbf{B}_k \widehat{\mathbf{v}}$, $\mathbf{w} = \mathbf{B}_k \widehat{\mathbf{w}}$. Since $(\mathbf{H}_k + \mu_k \mathbf{I})$ is block-diagonal, $(\mathbf{H}_k + \mu_k \mathbf{I})\widehat{\mathbf{v}} = \mathbf{u}$ can be decomposed into independent linear systems with coefficient matrices $H_s^k + \mu_k I$ for $s = 1, \ldots, K$ and we can proceed analogously for $(\mathbf{H}_k + \mu_k \mathbf{I})\widehat{\mathbf{w}} = \mathbf{g}^k$. Moreover, if we solve $(\mathbf{H}_k + \mu_k \mathbf{I})\widehat{\mathbf{v}} = \mathbf{u}$ by computing a factorization of $H_s^k + \mu_k I$, then the same factorization can be used to also solve the linear system for $\widehat{\mathbf{w}}$ and later to solve (15), so the computation of $\beta_k$ is not expensive.

Having $\beta_k$ computed as above, for the residual $\varphi_k(\beta_k)$ we have

$$\|\varphi_k(\beta_k)\|^2 = \|\mathbf{B}_k(\mathbf{H}_k + \mu_k \mathbf{I})^{-1}\mathbf{g}^k\|^2 + \beta_k^2 \|\mathbf{B}_k(I + (\mathbf{H}_k + \mu_k \mathbf{I})^{-1}\mathbf{B}_k)\mathbf{g}^k\|^2 \\ - 2\beta_k (\mathbf{g}^k)^\top (\mathbf{H}_k + \mu_k \mathbf{I})^{-\top} \mathbf{B}_k^\top \mathbf{B}_k (I + (\mathbf{H}_k + \mu_k \mathbf{I})^{-1}\mathbf{B}_k)\mathbf{g}^k. \tag{59}$$

If the vector $(\mathbf{H}_k + \mu_k \mathbf{I})^{-1}\mathbf{g}^k$ is in the null space of $\mathbf{B}_k$, we have that $\beta_k = 0$ and $|\varphi_k(\beta_k)\| = 0$, so in this case the direction $\mathbf{d}^k$ is equal to the Levenberg-Marquardt direction. If $\mathbf{g}^k$ is in the kernel of $\mathbf{B}_k$, then the residual $\|\varphi(\beta_k)\|$ is equal to $\|\mathbf{B}_k(\mathbf{H}_k + \mu_k\mathbf{I})^{-1}\mathbf{g}^k\|^2$ for any choice of the parameter $\beta$. If neither $(\mathbf{H}_k + \mu_k\mathbf{I})^{-1}\mathbf{g}^k$ nor $\mathbf{g}^k$ are in the null space of $\mathbf{B}_k$, then the optimal $\beta_k$ (58) is nonzero and so the right-hand side correction is effective in reducing the residual in the linear system. In general we have that

$$\|\varphi_k(\beta_k)\| \leq \|\mathbf{B}_k(\mathbf{H}_k + \mu_k\mathbf{I})^{-1}\mathbf{g}^k\| \leq \|\mathbf{B}_k\|\|(\mathbf{H}_k + \mu_k\mathbf{I})^{-1}\|\|\mathbf{g}^k\| \qquad (60)$$

so $\rho_k$ in (42) is bounded from above by $\|\mathbf{B}_k\|\|(\mathbf{H}_k + \mu_k\mathbf{I})^{-1}\|$. Taking into account Assumption 1, the definition of $\mu_k$ and the fact that this implies in particular $\mu_k \geq \frac{(b+1)^2}{1-b}\|\mathbf{J}_k\|^2$, we have

$$\rho_k \leq \|\mathbf{B}_k\|\|(\mathbf{H}_k + \mu_k\mathbf{I})^{-1}\| \leq \frac{M\|\mathbf{J}_k\|^2}{\mu_k} \leq \frac{M(1-b)}{(1+b)^2} < M. \qquad (61)$$

From the inequalities above, we have the dependence of the relative residual $\rho_k$ on the norm of the matrix $\|\mathbf{B}_k\|$, i.e., on the constant $M$ which measures the importance of the part that we disregard when approximating the Levenberg-Marquardt system with a block diagonal one. We can also notice that the residual is smaller for larger values of the damping parameter $\mu_k$.

# 5 Implementation and Numerical Results

In this section we present the results of a set of numerical experiments carried out to investigate the performance of the proposed method, compare it with classical Levenberg-Marquardt method and analyze the effectiveness of the right-hand side correction. For all the tests presented here we consider the case of Network Adjustment problems [17], briefly described in the Subsection 5.1. The LMS method is defined assuming that we can take advantage of the sparsity by suitable partition of variables and residuals and that we are able to apply the efficient right-hand-side correction as described in the Subsection 4.3, i.e., computing $\beta_k$ as in (58).

## 5.1 Least Squares Network Adjustment Problem

Consider a set of points $\{P_1, \ldots, P_n\}$ in $\mathbb{R}^2$ with unknown coordinates, and assume that a set of observations of geometrical quantities involving the points are available. Least Squares adjustments consists into using the available measurements to find accurate coordinates of the points, by minimizing the residual with respect to the given observations in the least squares sense.

We consider here network adjustment problems with three kinds of observations: point-point distance, angle formed by three points and point-line distance.

In order to be able to consider suitable increasing sizes, the problems are generated artificially, taking into account the information about average connectivity and structure of the network obtained from the analysis of real cadastral networks. The problems are generated as follows. Given the number of points $n$ we take $\{P_1, \ldots, P_n\}$ by uniformly sampling 25% of the points on a regular $2\sqrt{n} \times 2\sqrt{n}$ grid and we generate observations of the three kinds mentioned above until the average degree of the points is equal to 6. Each observation is generated by randomly selecting the points involved and generating a random number with Gaussian distribution with mean equal to the true measurement and given standard deviation. We use a standard deviation equal to 0.01 and 1 degree for distance and angle observations respectively. For all points we also add coordinates observations: for 1% of the points we use standard deviation 0.01, while for the remaining 99% we use standard deviation 1.

The optimization problem is defined as a weighted least squares problem

$$\min_{\mathbf{x} \in \mathbb{R}^N} \frac{1}{2} \sum_{j=1}^{m} r_j(\mathbf{x})^2 = \min_{\mathbf{x} \in \mathbb{R}^N} \frac{1}{2} \|\mathbf{R}(\mathbf{x})\|_2^2 \tag{62}$$

with $r_j(\mathbf{x}) = w_j^{-1} \widehat{r}_j(\mathbf{x})$, where $\widehat{r}_j$ is the residual function of the $j$-th observation and $w_j$ is the corresponding standard deviation.

In Figure 1 we present the spyplot of the matrix $\mathbf{J}^\top \mathbf{J}$ for a problem of size 35,000.

## 5.2 Comparison with Levenberg-Marquardt Method

In all the tests that follow we use a Python implementation of Algorithm LMS and classical LM method, and PyPardiso [9] to solve the sparse linear systems that arises at each iteration. All the tests were performed on a computer with Intel(R) Core(TM) i7-1165G7 processor @ 2.80GHz and 16.0 GB of RAM running Windows 10. All the methods that we consider have the same iteration structure. The main difference is the fact that while in LM method the linear system is solved directly using PyPardiso, in LMS we first perform the splitting and then use the same PyPardiso function to solve the resulting linear systems, therefore the comparisons in terms of time that we present are meaningful.

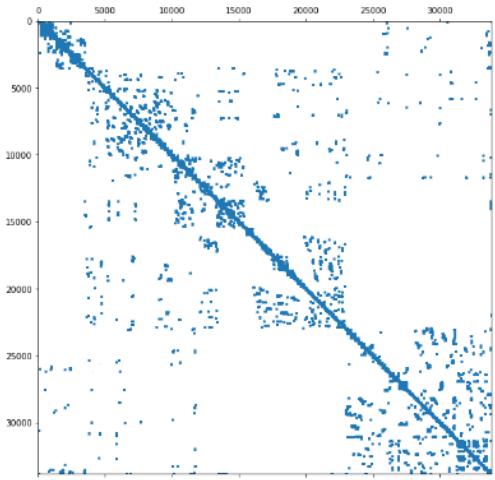

Figure 1: Sparsity plot of the coefficient matrix for N= 35,000

The partition of variables and residuals into sets $E_s, s = 1, \ldots, K$ is assumed to be given before application of LMS algorithm. To compute the partitioning of the variables, we use METIS [12] which, given a network and an integer $K > 1$ finds a partition of the vertices of the network into $K$ subsets of similar sizes, that approximately minimizes the number of edges between nodes in different subsets. The partition is computed by METIS in a multilevel fashion. Starting from a coarse representation of the graph, an initial partition is computed, projected onto a denser representation of the network and then refined. This process is repeated on a sequence of progressively more dense networks, up until the original graph. In all the tests that we considered, the time needed to compute the partitioning is negligible with respect to the overall computation time. This is in part due to the fact that the partitioning is computed only once at the beginning of the procedure and not repeated at each iteration.

We now consider a set of problems of increasing size and we solve each problem with the LMS method and correction coefficient $\beta_k$ computed as in (58). The problems are also solved with LM method. We consider problems with size between 20,000 and 120,000 and we plot the time taken by the two methods to reach termination. Both methods use as initial guess the coordinate observations available in the problem description and they stop when at least $68\%, 95\%, 99.5\%$ of the residuals is smaller than 1, 2 and 3 times the standard deviation respectively. The obtained results are in the first plot of Figure 2. To give a better comparison, in the second plot we extend the size of the problems solved with the proposed method up to 1 million variables. Clearly, LM method could not cope with such large problems (in our testing environment) while LMS successfully solved problems of increasing dimensions up to final value of 1 million variables. In Figure 3 we have the log-log plot of the time necessary to solve each problem, compared with different rates of growth. For the method with $K > 1$, a small number of values of the parameter $K$ was tested and the best one was selected to perform the comparison. The value $K$ used at each dimension is reported in Figure 8.

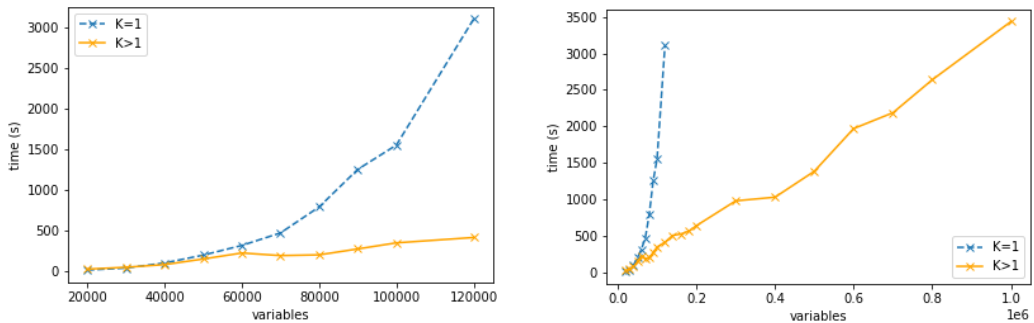

Figure 2: Time comparison between classical Levenberg-Marquardt ($K = 1$) and the proposed method ($K > 1$) with optimal $\beta_k$ for problems of increasing size.

From Figure 2 one can see that the LMS method with $K > 1$ results in a significant reduction of the time necessary to reach the desired accuracy, compared to Levenberg-Marquardt method. Moreover, from the second plot of Figure 2 and from Figure 3 we can notice that, on the problems that we considered, the time taken by the proposed method grows approximately as $n^{1.3}$, which suggests the fact that the method discussed in this paper would be suitable for problems of very large dimensions.

To better understand the behaviour of the method, in Figure 4 we plot the mentioned percentages and the value of the relative residual $F_k/F_0$ at each iteration, for a problem of size $N = 10^5$ and $K = 15$. For the same problem, in Figure 5 we plot the distribution plot of the coordinate error with respect to the true solution, for the initial guess and the estimated solution (left and right-hand plot, respectively).

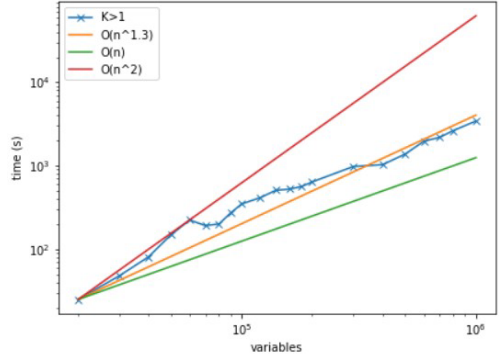

Figure 3: Dependence of time on the size of the problem - loglog plot, $N \leq 10^6$.

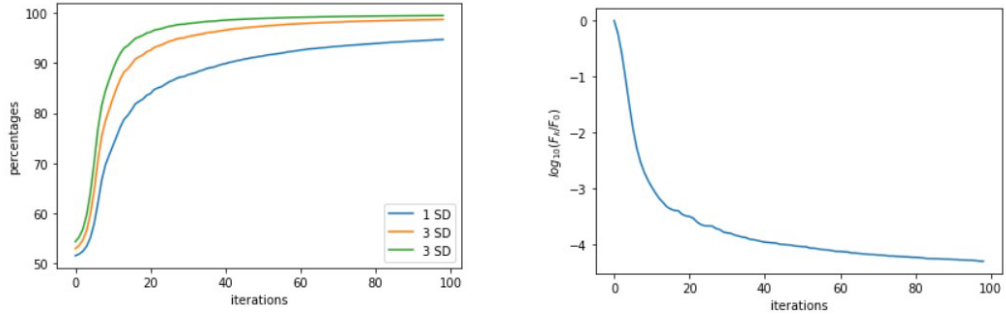

Figure 4: The percentage of points within standard deviation measure (left) and the relative residual (right) per iteration for $N = 10^5$ and $K = 15$.

## 5.3 Influence of the parameters $K$ and $\beta_k$

Let us now study how the number of subproblems $K$ influences the performance of the method. We consider two problems with 100,000 and 200,000 variables respectively and we solve them with the proposed algorithm for a set of increasing values of $K$. For each considered $K$ we plot in Figure 6 the time taken by the method to arrive at the desired accuracy. The initial guess and the stopping criterion are defined as in the previous test.

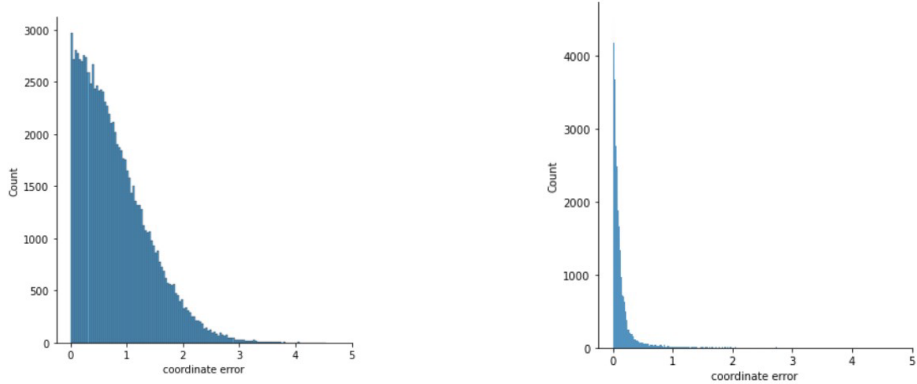

Figure 5: Distribution plot of the coordinate error for all points, at the initial guess (left) and at the final iterate (right).

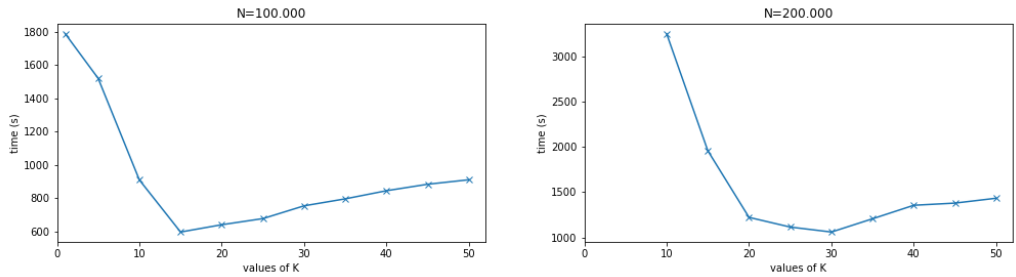

Figure 6: Time to compute a solution for optimal $\beta_k$ and different values of $K$, $N = 10^5$ and $N = 2 \cdot 10^5$

One can notice that the time decreases as $K$ increases up to an optimal value ($K = 15$ for the first problem and $K = 30$ for the second one) after which the time starts to increase again. The reason behind this behavior is that larger values of $K$ yield smaller linear system and therefore cheaper iterations, but also less accurate search direction $\mathbf{d}^k$ resulting in a larger number of iterations necessary to achieve the desired accuracy. For large values of $K$ the increase in the number of iterations outweights the saving in the solution of the linear system and the overall computation cost increases. Finally we can notice that, despite the existence of an optimal value of the parameter $K$, it appears from this test that there exists an interval of values for which the cost of the method is comparable. This suggests that fine-tuning of the parameter $K$ is not necessary and that, given a problem, choosing $K$ according to the number of variables should be enough to achieve good performance.

To see that the proposed right-hand side correction improves the performance of the method, we repeat the test presented in Subsection 5.2 for $N = 10^6$, but the comparison is here carried out with the case $\beta_k = 0$ that is, when the linear system is approximated as in (13) but no right-hand side correction is applied.

For both methods, a few different values of the parameter $K$ were tested. In Figure 7 we report the time needed by the two methods to satisfy the convergence criterion, for the best $K$ among the considered values. In figure 8 we plot for each method and each size the value of the parameter $K$ corresponding to the timings in Figure 7.

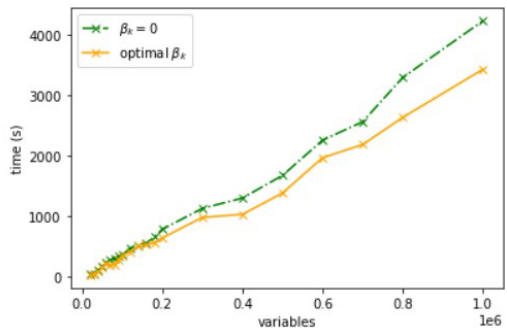

Figure 7: Time to compute the solution with $\beta_k = 0$ and optimal $\beta_k$.

We can see that applying the proposed right-hand side correction effectively reduce the time necessary to satisfy the stopping condition. From Figure 7 one can notice that the optimal $K$ for the method with right-hand side correction is generally higher than the method without correction. These two results together suggests that the method with right-hand side correction is able to achieve a better performance because it allows the set of variables to be partitioned into smaller subsets, which implies a faster computation of the direction at each iteration, before incurring into a decrease in the performance due to the additional number of iterations necessary to reach the desired accuracy.

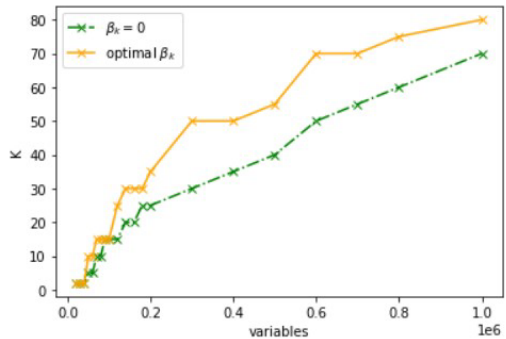

Figure 8: Selected values of the number of subproblems $K$

# 6 Conclusions

We presented a method of inexact Levenberg-Marquardt type for sparse problems of very large dimension. Assuming that the problem is nearly separable, i.e., sufficiently sparse such that each component of the residual depends only on few variables, the proposed methods is defined through splitting into a set of $K$ independent systems of equations of smaller dimension. The decoupling is done taking into account dense diagonal blocks of LM system and disregarding hopefully very sparse off-diagonal blocks. To compensate for disregarded off-diagonal block we introduced a correction on the right-hand side of the system in such way that decoupling is maintained but information contained in the off-diagonal matrix is preserved in a computationally affordable way, using a single parameter that can be computed in the same fashion - solving a sequence of small dimensional systems of linear equations. The key idea is that solving $K$ systems of smaller dimensions, that can be done sequentially or in parallel, is significantly cheaper than solving a large system of linear equations even if the system is sparse.

The presented algorithm is globally convergent under the set of standard assumptions for a suitable choice of regularization parameter in LM system. In fact the global convergence does not rely on separability assumption at all as one can show that the direction computed by decoupled sequence of LM systems is descent direction. To achieve global convergence we rely on linesearch and regularization parameter update by a trust-region like scheme, similarly to [11]. Local linear convergence is proved under the standard conditions and assuming that the residual of linear system is small enough in each iteration. Hence, the near-separability assumption plays a role in local convergence. To achieve small residuals for the decoupled problem we rely heavily on the right-hand side correction and discuss the optimal choice of parameter that is employed in the correction. Theoretical considerations are supported by numerical examples. We consider the network adjustment problem on simulated data, inspired by a real-world problem of cadaster maps, of growing size and with the proposed method solve problems of up to one million of variables. Comparison with the classical LM is presented and it is shown that the proposed method is significantly faster and able to cope with large dimensions. The experiments reported in this paper are done in sequential way while the parallel implementation will be a subject of further research.

# Funding

This work is supported by the European Union's Horizon 2020 programme under the Marie Skłodowska-Curie Grant Agreement no. 812912. The work of Krejić is partially supported by the Serbian Ministry of Education, Science and Technological Development.